\documentclass[12pt]{article} 
\usepackage[a4paper,left=2.54cm, right=2.54cm, bottom=2.54cm, top=2.54cm]{geometry}
\usepackage{amssymb,float,latexsym,caption,subcaption}
\usepackage{amsmath,amsthm}
\usepackage{booktabs}
\usepackage{commath}
\usepackage[utf8]{inputenc}
\usepackage[english]{babel}
\usepackage{graphicx}
\usepackage{epstopdf,epsfig}
\usepackage[labelfont={bf}]{caption}
\usepackage{color}
\usepackage[T1]{fontenc}
\usepackage[mathscr]{eucal}
\usepackage{authblk}
\usepackage[round,sort]{natbib}
\usepackage{titlesec}
\newcommand{\secfnt}{\fontsize{14}{14}}
\newcommand{\ssecfnt}{\fontsize{14}{14}}
\titleformat{\chapter}[display]
{\normalfont\chapfnt\bfseries}{\chaptertitlename\ \thechapter}{12pt}{\chapfnt}
\titleformat{\section}
{\normalfont\secfnt\bfseries}{\thesection}{1em}{}
\titleformat{\subsection}
{\normalfont\ssecfnt\bfseries}{\thesubsection}{1em}{}
\titlespacing*{\chapter} {0pt}{50pt}{40pt}
\titlespacing*{\section} {0pt}{3.5ex plus 1ex minus .2ex}{2.3ex plus .2ex}
\titlespacing*{\subsection} {0pt}{3.25ex plus 1ex minus .2ex}{1.5ex plus .2ex}
\newtheorem{theorem}{Theorem}[section]
\newtheorem{example}{Example}[section]
\newtheorem{remark}{Remark}[section]
\newtheorem{definition}{Definition}[section]
\newtheorem{assumption}{Assumption}[section]

	\author{Tareq Alodat$^1$ and Andriy Olenko$^{1}\footnote{Contact: A.Olenko email \href{mailto:a.olenko@latrobe.edu.au}{a.olenko@latrobe.edu.au}}$}
\affil{\small$^{1}$Department of Mathematics and Statistics, La Trobe University, Melbourne, VIC, 3086, Australia} 
\date{}
\title{\bf Asymptotic Behaviour of Discretised Functionals of Long-Range Dependent Functional Data}
\numberwithin{equation}{section}

\begin{document}
\bibliographystyle{plainnat}
\maketitle
{The paper studies the asymptotic behaviour of weighted functionals of long-range dependent data over increasing observation windows. Various important statistics, including sample means, high order moments, occupation measures can be given by these functionals. It is shown that in the discrete sampling case additive functionals have the same asymptotic distribution as the corresponding integral functionals for the continuous functional data case. These results are applied to obtain non-central limit theorems for weighted additive functionals of random fields. As the majority of known results concern the discrete sampling case the developed methodology helps in translating these results to functional data without deriving them again. Numerical studies suggest that the theoretical findings are valid for wider classes of long-range dependent data.}

\begin{flushleft}
\end{flushleft}
\section{Introduction}\label{sec1}
Recent advances in technology allowed collecting big data at high frequency (effectively continuous) rates that led to the ubiquity of functional data (samples of curves or surfaces)~(\cite{ramsay2005,wang2016functional}). Handling such new complex data is essential in various applications, for example, in earth, environmental, ecological sciences, cosmology and image analysis. However, most of classical statistical models and results were developed for discretely sampled~data.

Discretisation and corresponding additive models are often used as powerful dimension reduction tools in the analysis of functional data which are intrinsically infinite dimensional. The discretisation  is a common strategy for approximating statistics of such data (see, for example, $\S\ 6.4.3$ in~\cite{ramsay2005}). Also, in practice, the functional curves or surfaces are often observed only at a finite number of points. 

Note, that various statistics of functional data can be expressed by integral functionals of these data or their transformations. For instance, some well-known examples of such statistics include sample moments and sample sojourn measures (Minkowski functionals)~(see \cite{leonenko2014sojourn}). Another popular model in various applications (especially in engineering and signal processing) is stochastic processes that are obtained as outputs of filters, i.e. defined mathematically by a convolution integral operator.
Another example is functional linear regression models defined by weighted integral functionals. These models found numerous statistical applications in medicine, linguistics, chemometrics~(see~\cite{ramsay2005};~\cite{crambes2009smoothing};~\cite{zhanganalysis}).

In all above applications it is usually assumed that the discretisation error is negligible with respect to the estimation error. However, there are almost no known results that rigorously prove it. This paper addresses this problem and investigates discretisation errors for weighted functionals of long-range dependent spatial processes. Their rates of decay for the case of increasing observation windows are found. It is shown that both additive and integral functionals converge to the same limit distribution. It is proved that these distributions are non-Gaussian. These results provide a constructive method for determining the number of discretisation nodes for a given accuracy.

Various results in statistical inference of random fields were first obtained by~\cite{yadrenko1983spectral}. Recently, considerable attention has been paid to asymptotic behaviour of non-linear statistics of random processes and fields~(see~\cite{ivanov2008semiparametric};~\cite{bai2013multivariate};~\cite{leonenko2014sojourn};~\cite{anh2017rate} and the references therein). Direct probability techniques were used to study these statistics, for example, in regression models. Asymptotic distributions of these statistics were discussed by~\cite{ivanov1989statistical} and it was shown that central and non-central limit theorems hold for particular models. However, no results about discretisations were given.

There are many practical situations in which non-Gaussian random processes and fields are appropriate for statistical modeling. We deal with an important class of models defined by non-linear functions of Gaussian random fields. This class is widely used in modeling non-Gaussian data. It can be analysed using Wiener chaos expansions that give good data approximations in many cases (see~\cite{de1997bayesian,vio2001numerical}).

This research deals with asymptotic behaviour of integral and additive non-linear functionals of random fields with long-range dependence. Long-range dependence is an empirical phenomenon which has been observed in different applied fields including cosmology, economics, geophysics, air pollution, image analysis, earth sciences, just to mention a few examples. For this reason, great effort has been devoted to studying models based on long-range dependent random fields 
(see~\cite{ivanov1989statistical};~\cite{wackernagel1998multivariate};~\cite{doukhan2002theory};~\cite{frias2008}). Weighted functionals of long-range dependent random fields were considered in~\cite{olenko2013limit,ivanov2013limit,ivanov1989statistical}.
These functionals can have non-Gaussian limits that are known as Hermite or Hermite-Rosenblatt distributions~(\cite{rosenblatt1961independence};~\cite{taqqu1975weak};~\cite{dobrushin1979non,taqqu1979convergence}). Their asymptotic distributions can be characterised by either multiple Wiener-It{\^o} integrals representations or characteristic functions (see~\cite{taqqu1979convergence,dobrushin1979non};~\cite{leonenko2006weak}).

In various applications, it is natural to consider statistics of random fields and to study their limit behaviour over increasing spatial windows. In these cases, integrals of non-linear functionals of spatial functional data and additive non-linear functionals for discrete observations on a bounded region were studied in numerous papers (see, for example,~\cite{major1981multiple};~\cite{leonenko2014sojourn};~\cite{anh2015rate,anh2017rate}). When we deal with the asymptotic behaviour of discretised functionals of functional data, it is important to know how asymptotics of these integrals are related to additive functionals. To the best of our knowledge only particular cases of this correspondence have been addressed in~\cite{leonenko2006weak} and~\cite{weak2017alodat} for rectangular observation windows. However, in many applications spatial data is not necessarily available over rectangles, but rather over irregularly-shaped regions~(\cite{cressie1993statistics,lahiri1999prediction}). Therefore, it is important to obtain theoretical results about asymptotics for general types of observation windows. In this paper we extend results of~\cite{leonenko2006weak} and~\cite{weak2017alodat} under more general conditions. More precisely, we consider weighted functionals of random fields of the form 
	\begin{align*}
d_{r}^{-1}\int_{\Delta_{n}(r)}g(x)H_{\kappa}(\xi(x))dx,\quad r\to\infty,
\end{align*}
where $\xi(x),\ x\in\mathbb{R}^{n}$, is a long-range dependent random field, $g(x)$ and $H_{\kappa}(\cdot)$ are non-random functions, $\Delta_{n}\subset\mathbb{R}^{n}$ is an observation window and $d_{r}^{-1}$ is a normalising factor. We show that these integrals and the corresponding discretised versions have same non-Gaussian limit distributions. 

The article is organised as follows. In Section~\ref{sec2} we introduce main notations and outline necessary background from the theory of random fields. In Section~\ref{sec3} we recall some assumptions and auxiliary results from the spectral and correlation theory of random fields. In Section~\ref{sec4} we study the case of two-dimensional functionals. Section~\ref{sec5} gives a general  multidimensional version of the results. Proofs are provided in Section~\ref{secproofs}. In Section~\ref{sec6} some simulations studies are presented to confirm theoretical results. Conclusions and directions for future research are presented in Section~\ref{sec7}.
\section{Definitions and Notations}\label{sec2}
	This section provides basic definitions and notations that are used in this article.
	
	In what follows $|\cdot|$, $\Vert \cdot \Vert$, $\lfloor\cdot\rfloor$ and $\lceil\cdot\rceil$ are used for the Lebesgue measure, the Euclidean distance in $\mathbb{R}^{n}$, the floor and ceiling functions, respectively. The symbols $C$, $\varepsilon$ and $\delta$ (with subscripts) will be used to denote constants that are not important for our discussion. Note, that the same symbol~$C$ may be used for different constants appearing in the same proof. For a set $A\subset\mathbb{R}^{n},\ n\geq1$, we denote by $A^{\circ}$ and $A^{c}$ the interior and the exterior of the set $A$ respectively. Moreover, it is assumed that all random variables are defined on a fixed probability space $\left(\Omega,\mathfrak{F},\mathbb{P}\right)$.
	
	We consider a measurable mean-square continuous zero-mean homogeneous isotropic real-valued random field $\xi\left(x\right),\ x\in\mathbb{R}^{n}$, with the covariance function
	\begin{equation*}
	B\left(r\right)=B(\Vert x\Vert):=\mathbb{E}\left(\xi(0)\xi(x)\right),\quad x\in\mathbb{R}^{n}.
	\end{equation*}
	
	It is well known that there exists a bounded non-decreasing function $\Phi\left(u\right),\ u\geqslant 0$, (see~\cite{yadrenko1983spectral};~\cite{ivanov1989statistical}) such that
	\begin{equation*}
	B\left(r\right)=\big(2/r\big)^{(n-2)/2}\Gamma\big(n/2\big)\int_{0}^{\infty}J_{(n-2)/2}(ru)u^{(2-n)/2}d\Phi\left(u\right),
	\end{equation*}
	where $J_{v}(\cdot)$ is the Bessel function of the first kind of order $v>-1/2$. 
	
	The function $\Phi\left(\cdot\right)$ is called the isotropic spectral measure of the random field $\xi\left(x\right),$ $x\in\mathbb{R}^{n}$. If there exists a function $\varphi(u),\ u\in[0,\infty)$, such that
	\begin{equation*}
	u^{n-1}\varphi(u)\in L_{1}([0,\infty)),\quad\Phi(u)=\dfrac{2\pi^{n/2}}{\Gamma(n/2)}\int_{0}^{u}z^{n-1}\varphi(z)dz,
	\end{equation*}
	then the function $\varphi(\cdot)$ is called the isotropic spectral density of the field $\xi\left(x\right)$.
	
	The field $\xi\left(x\right)$ with an absolutely continuous spectrum has the following isonormal spectral representation
	\begin{align}\label{eq0}
	\xi\left(x\right)=\int_{\mathbb{R}^{n}}e^{i\langle\lambda,x \rangle}\sqrt{\varphi(\Vert\lambda\Vert)}W(d\lambda),
	\end{align}
	where $W(\cdot)$ is the complex Gaussian white noise random measure on $\mathbb{R}^{n}$~(see~\cite{yadrenko1983spectral}; \cite{ivanov1989statistical}).

The Hermite polynomials $H_{m}(x),\ m\geq 0$, are defined by	
	\begin{equation*}
	H_{m}(x):=(-1)^{m}\exp\left(\dfrac{x^{2}}{2}\right)\dfrac{d^{m}}{dx^{m}}\exp\left(-\dfrac{x^{2}}{2}\right).
	\end{equation*}
	The first few Hermite polynomials are
	$
	H_{0}(x)=1, H_{1}(x)=x,  H_{2}(x)=x^{2}-1, H_{3}(x)=x^{3}-3x.
	$
	
	The Hermite polynomials $H_{m}(x),\ m\geq 0$, form a complete orthogonal system in the Hilbert space $L_{2}\left(\mathbb{R},\phi(\omega)d\omega\right):=\{G: \int_{\mathbb{R}}G^{2}(\omega)\phi(\omega)d\omega<\infty\}$, where 
	$\phi(\omega)$ is the probability density function of the standard normal distribution. An arbitrary function $G(\omega)\in L_{2}\left(\mathbb{R},\phi(\omega)d\omega\right)$ admits the mean-square convergent expansion
	\begin{equation}\label{eq2}
	G(\omega)=\sum_{j=0}^{\infty}\dfrac{C_{j}H_{j}(\omega)}{j!},\quad C_{j}:=\int_{\mathbb{R}}G(\omega)H_{j}(\omega)\phi(\omega)d\omega. 
	\end{equation}
	By Parseval's identity it holds 
	\begin{equation*}
	\sum_{j=0}^{\infty}\dfrac{C_{j}^{2}}{j!}=\int_{\mathbb{R}}G^{2}(\omega)\phi(\omega)d\omega.
	\end{equation*}
	\begin{definition}\label{def1}\rm(\cite{taqqu1975weak}) Let $G(\cdot)\in L_{2}\left(\mathbb{R},\phi(\omega)d\omega\right)$. Assume that there exists an integer $\kappa\geqslant 1$, such that $C_{j}=0$ for all $0<j\leq \kappa-1$, but $C_{\kappa}\neq 0$. Then $\kappa$ is called the Hermite rank of $G(\cdot)$ and is denoted by $H rankG(\cdot).$
	\end{definition}
		Note, that by (2.1.8) in~\cite{ivanov1989statistical} we get $\mathbb{E}\left(H_{m}(\xi(x))\right)=0$ and
		\begin{equation}\label{eq1}
		\mathbb{E}\left(H_{m_{1}}(\xi(x))H_{m_{2}}(\xi(y))\right)=\delta_{m_{1}}^{m_{2}}m_{1}!B^{m_{1}}(\Vert x-y\Vert),\quad x,y\in\mathbb{R}^{n},
		\end{equation}
		where $\delta_{m_{1}}^{m_{2}}$ is the Kronecker delta function.
	\begin{definition}\label{def3}\rm(\cite{bingham1989regular}) A measurable function $L:(0,\infty)\rightarrow (0,\infty)$ is slowly varying at infinity if for all $t>0$,
$\lim_{r\to\infty} L(tr)/L(r)=1.$
	\end{definition}
\section{Assumptions and Auxiliary Results}\label{sec3}
This section gives some assumptions and results from the spectral and correlation theory of random fields that will be used in the following sections.
	\begin{assumption}\label{ass1}
		Let $\xi(x),\ x\in\mathbb{R}^{n}$, be a homogeneous isotropic Gaussian random field with $\mathbb{E}\xi(x)=0$ and the covariance function $B(x)$, such that $B(0)=1$ and
		$$B(x)=\mathbb{E}\left(\xi\left(0\right)\xi\left(x\right)\right)=\Vert x\Vert^{-\alpha}L_{0}\left(\Vert x\Vert\right),\quad \alpha>0,$$\
		where $L_{0}\left(\Vert\cdot\Vert\right)$ is a function slowly varying at infinity.  
	\end{assumption}
	
	If $\alpha\in\left(0,n\right)$, then the covariance function $B(x)$ satisfying Assumption \ref{ass1} is not integrable, which corresponds to the long-range dependence case~(\cite{anh2015rate}).
	
	The notation $\Delta_{n}\subset\mathbb{R}^{n}$ will be used to denote a Jordan-measurable compact bounded set, such that $|\Delta_{n}|>0$, and $\Delta_{n}$ contains the origin in its interior. Let $\Delta_{n}(r),\ r>0$, be the homothetic image of the set $\Delta_{n}$, with the centre of homothety at the origin and the coefficient $r>0$, that is $\vert\Delta_{n}(r)\vert=r^{n}\vert\Delta_{n}\vert$ and $\Delta_{n}=\Delta_{n}(1)$. For $c\in\mathbb{R}$ and $v\in\mathbb{R}^{n},\ n\geq 1$, we define $c\Delta_{n}:=\set{cx:x\in\Delta_{n}}$ and $\Delta_{n}-v:=\set{x-v:x\in\Delta_{n}}$. 
	
	Let $HrankG(\cdot)=\kappa$. Denote the random variables $K_{r}$ and $K_{r,\kappa}$ by
	$$K_{r}:=\int_{\Delta_{n}(r)}G\left(\xi\left(x\right)\right)dx\quad \text{and}\quad K_{r,\kappa}:=\dfrac{C_{\kappa}}{\kappa!}\int_{\Delta_{n}(r)}H_{\kappa}\left(\xi\left(x\right)\right)dx,$$
	where $C_{\kappa}$ is given by (\ref{eq2}).
	\begin{theorem}\label{theo1}{\rm(\cite{leonenko2014sojourn})}
		Suppose that $\xi\left(x\right),\ x\in\mathbb{R}^{n}$, satisfies Assumption~{\rm\ref{ass1}} and $H rankG(\cdot)=\kappa\geq 1$. If a limit distribution exists for at least one of the random variables
		$$\dfrac{K_{r}}{\sqrt{Var K_{r}}}\quad and\quad \dfrac{K_{r,\kappa}}{\sqrt{Var K_{r,\kappa}}},$$
		then the limit distribution of the other random variable also exists, and the limit distributions coincide when $r\to\infty$. 
	\end{theorem}
	
	By Theorem \ref{theo1} it is enough to study $K_{r,\kappa}$ to get asymptotic distributions of $K_{r}$. Therefore, we restrict our attention only to $K_{r,\kappa}$.
	
	\begin{assumption}\label{ass2}
		The random field $\xi\left(x\right),\  x\in\mathbb{R}^{n}$, has the isotropic spectral density
		\begin{equation*}
		\varphi\left(\Vert \lambda\Vert\right):=c_{1}\left(n,\alpha\right)\Vert\lambda\Vert^{\alpha-n}L\left(\dfrac{1}{\Vert  \lambda\Vert}\right),
		\end{equation*}
		where $\alpha\in(0,n),\ c_{1}\left(n,\alpha\right):=\Gamma\left((n-\alpha)/{2}\right)/2^{\alpha}\pi^{n/2}\Gamma\left(\alpha/2\right),$ and $L(\Vert\cdot\Vert)\sim L_{0}(\Vert\cdot\Vert)$ is a locally bounded function which is slowly varying at infinity.
	\end{assumption}
	
	One can find more details on relations between Assumptions~\ref{ass1} and~\ref{ass2} in~\cite{anh2017rate}.

	The function $K_{\Delta_{n}}\left(x\right)$ will be used to denote the Fourier transform of the indicator function of the set $\Delta_{n}$, i.e.
	\begin{align*}
	K_{\Delta_{n}}\left(x\right):=\int_{\Delta_{n}}e^{i\langle u,x\rangle}du,\quad x\in\mathbb{R}^{n}.
\end{align*}
\begin{theorem}\label{theo2}{\rm(\cite{leonenko2014sojourn})}
	Let $\xi\left(x\right),\ x\in\mathbb{R}^{n}$, be a homogeneous isotropic Gaussian random field. If Assumptions {\rm\ref{ass1}} and {\rm\ref{ass2}} hold, $\alpha\in\left(0,n/\kappa\right)$, then for $r \rightarrow\infty$ the random variables 
		$$X_{r,\kappa}(\Delta_{n}):=r^{\kappa\alpha/2-n}L^{-\kappa/2}(r)\int_{\Delta_{n}(r)}H_{\kappa}\left(\xi(x)\right)dx$$\
		converge weakly to  
		\begin{align*}
		X_{\kappa}(\Delta_{n}):=c_{1}^{\kappa/2}(n,\alpha)\int_{\mathbb{R}^{n\kappa}}^{\prime}K_{\Delta_{n}}\left(\lambda_{1}+\cdots+\lambda_{\kappa}\right)\dfrac{W(d\lambda_{1})\cdots W(d\lambda_{\kappa})}{\Vert\lambda_{1}\Vert^{(n-\alpha)/2}\cdots\Vert\lambda_{\kappa}\Vert^{(n-\alpha)/2}}.
		\end{align*}
		Here $\int_{\mathbb{R}^{n\kappa}}^{\prime}$ denotes the multiple Wiener-It{\^o} integral with respect to a Gaussian white noise measure, where the diagonal hyperplanes $\lambda_{i}=\pm\lambda_{j},\ i,j=1,\dots,\kappa,\ i\neq j$, are excluded from the domain of integration.
	\end{theorem}

Below we present a limit theorem and the corresponding assumptions on the weight function in the integral functionals from~\cite{ivanov1989statistical}. These results will be generalised in the subsequent sections. The obtained results on asymptotic equivalence of additive and integral functionals of random fields will be used to obtain limit theorems for the case of discrete observations. 	
\begin{assumption}\label{ass2weighted}{\rm(\cite{ivanov1989statistical})} 
	Let $\vartheta(x)=\vartheta(\Vert x\Vert)$ be a radial continuous function that is positive for $\Vert x\Vert>0$ and such that for $\alpha\in\left(0,n/\kappa\right)$
	\begin{align*}
   \lim_{r\to\infty}\int_{\Delta_{n}}\int_{\Delta_{n}}\dfrac{\vartheta(r\Vert x\Vert)\vartheta(r\Vert y\Vert)dxdy}{\vartheta^{2}(r)\Vert x-y\Vert^{\alpha\kappa}}\in\left(0,\infty\right).
	\end{align*}
	\end{assumption}

Let $u(\Vert\lambda\Vert):=c_{1}\left(n,\alpha\right)L\left(\frac{1}{\Vert \lambda\Vert}\right)$, where $L(\cdot)$ is from Assumption \ref{ass2}. In Section~2.10 in~\cite{ivanov1989statistical} the case when the function $u(\Vert\lambda\Vert)$ is continuous in a neighborhood of zero, bounded on $ (0,\infty) $  and $u(0)\neq0$, was studied. It was assumed that there is a function $\bar{\vartheta}(\Vert x\Vert)$ such that
\begin{align*}
\int_{\mathbb{R}^{n\kappa}}\prod_{j=1}^{\kappa}\Vert \lambda_{j}\Vert^{\alpha-n}\abs{\int_{\Delta_{n}}e^{i\langle \lambda_{1}+\dots+\lambda_{\kappa},x \rangle}\bar{\vartheta}(x)dx}^{2}\prod_{j=1}^{\kappa}d\lambda_{j}<\infty
\end{align*}
and
\[
\lim_{r\to\infty}\int_{\mathbb{R}^{n\kappa}}\bigg|\int_{\Delta_{n}}e^{i\langle \lambda_{1}+\dots+\lambda_{\kappa},x \rangle}
\bigg(\frac{\vartheta(r\Vert x \Vert)}{\vartheta(r)}\prod_{j=1}^{\kappa}\sqrt{\frac{u(\Vert\lambda_{j}\Vert r^{-1})}{u(0)}}-\bar{\vartheta}(x)\bigg)dx\bigg|^{2}\prod_{j=1}^{\kappa}\Vert \lambda_{j}\Vert^{\alpha-n}\prod_{j=1}^{\kappa}d\lambda_{j}=0.
\]

Under these assumptions the following result was obtained.
\begin{theorem}\label{theo3}{\rm(\cite{ivanov1989statistical})}
	If Assumption {\rm\ref{ass2weighted}} holds, then for $r \rightarrow\infty$ the random variables
	$$Y_{r,\kappa}:=\frac{1}{r^{n-\kappa\alpha/2}\vartheta(r)u^{\kappa/2}(0)}\int_{\Delta_{n}(r)}\vartheta(\Vert x\Vert)H_{\kappa}\left(\xi(x)\right)dx$$
	converge weakly to
	\begin{align*}
	Y_{\kappa}:=\int_{\mathbb{R}^{n\kappa}}^{\prime}K_{\Delta_{n}}\left(\lambda_{1}+\cdots+\lambda_{\kappa};\bar{\vartheta}\right)\dfrac{\prod_{j=1}^{\kappa} W(d\lambda_{j})}{\prod_{j=1}^{\kappa}\Vert\lambda_{j}\Vert^{(n-\alpha)/2}},
	\end{align*}  
where
$\alpha\in\big(0,\min\left(\frac{n}{\kappa},\frac{n+1}{2}\right)\big)$
and
$
K_{\Delta_{n}}\left(\lambda;\bar{\vartheta}\right):=\int_{\Delta_{n}}e^{i\langle \lambda,x \rangle}\bar{\vartheta}(x)dx.
$
\end{theorem}
 Note, that the above result is a specification of the results on convergence to stochastic processes in~\cite{ivanov1989statistical} where functionals over the observation windows $\Delta_{n}(rt^{1/n}),\ t\in[0,1]$, were studied. For simplicity, this paper deals with a particular case when $t=1$. But the results of this paper can be easily extended to the case of general $\Delta_{n}(rt^{1/n})$ and convergence to stochastic processes on $t\in[0,1]$.
 \section{Main Results}
 \subsection{Two-dimensional Case}\label{sec4}
In this section we consider integrals of two-dimensional random fields over a bounded increasing observation window $\Delta_{2}(r)\subset\mathbb{R}^{2},\ r>0$. We show that the limit distributions of these integrals and their corresponding additive functionals coincide. 

\begin{figure}[h]
	\centering
	\begin{subfigure}[b]{0.4\textwidth}
		\includegraphics[width=1\textwidth,height=5.1cm,
		trim={-1cm 0 0 0},clip]{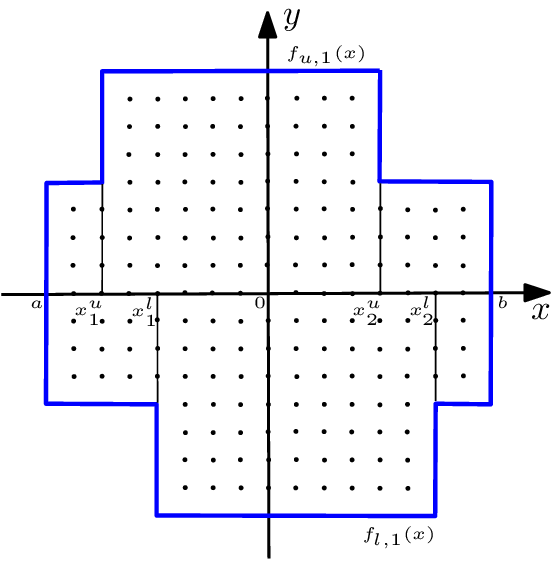}  \vspace{-0.5cm}
		\caption{}
		\label{fig1:case2a}
	\end{subfigure}
	\hfill  
	\begin{subfigure}[b]{0.4\textwidth}
		\includegraphics[width=0.97\textwidth, height=5.cm,trim={.1cm 0 0 0},clip]{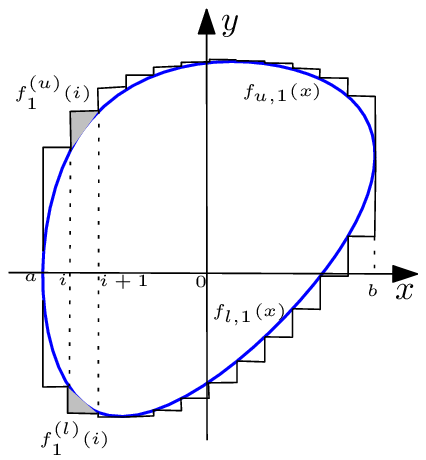} \vspace{0.5cm}
		\caption{}
		\label{fig1:case2b}
	\end{subfigure} \hspace{2cm}
	\vspace{-.3cm}
	\caption{(a) Two-dimensional set $\Delta_{2}$ with a non-smooth boundary, (b) Two-dimensional set $\Delta_{2}$ and its $f_{1}^{(l)}(x)$ and $f_{1}^{(u)}(x)$. The shaded areas are $S_{q,1}(i)\cap\Delta_{2}^{c}$.}\label{fig1}
\end{figure}

Our setup is as follows. We assume that  the set $\Delta_{2}$ can be represented as
$$\Delta_{2}=\{(x,y)\in\mathbb{R}^{2}:a\leq x\leq b,\  f_{l,1}(x)\leq y\leq f_{u,1}(x)\},$$ where $a=\displaystyle\min_{(x,y)\in\Delta_{2}}x,\  b=\displaystyle\max_{(x,y)\in\Delta_{2}}x$, $f_{l,1}(x)<f_{u,1}(x),\ x\in(a,b)$, and $f_{q,1}(x),\ q\in\{l,u\}$, are smooth functions (i.e. $f_{q,1}(x)\in\mathcal{C}^{1},\ q\in\{l,u\}$, where $\mathcal{C}^{1}$ is a class of functions with continuous first derivatives) except of the sets $ M_{q}=\{x_{1}^{q},\dots,x_{k_{q}}^{q}\}\subset[a,b]$, where these functions have finite jumps. Here $k_{q}$ is the number of jumps of $f_{q,1}(\cdot)$, see, for example, Figure~\ref{fig1:case2a}. That is $x^{q}_{j}\in M_{q}$, $j=1,\dots,k_{q}$, are jump points of the functions $f_{q,1}(\cdot),\ q\in\{l,u\},$ if $\lim\limits_{x\to x^{q+}_{j}}f_{q,1}(x)\neq \lim\limits_{x\to x^{q-}_{j}}f_{q,1}(x)$ but $\lim\limits_{x\to x^{q+}_{j}}f_{q,1}(x)$ and $\lim\limits_{x\to x^{q-}_{j}}f_{q,1}(x)$ both exist. Note that, by the  homothety of $\Delta_{2}(r),\ r>0$, the set $\Delta_{2}(r)$ can be represented as 
$\Delta_{2}(r)=\{(x,y)\in\mathbb{R}^{2}:ar\leq x\leq br, f_{l,r}(x)\leq y\leq f_{u,r}(x)\},$
where $f_{q,r}(x)=rf_{q,1}(x/r),\ q\in\{l,u\}$. As $\Delta_{2}$ contains the origin in its interior, it follows that $ar<0< br$, $\lfloor ar\rfloor\to-\infty$, and $\lceil br\rceil\to\infty,$ as $r\to\infty$.
\raggedbottom
Let $\xi(x,y),\ x, y\in\mathbb{R}$, be a real-valued homogeneous isotropic Gaussian random field satisfying Assumptions~\ref{ass1} and \ref{ass1}. We investigate the integrals
	\begin{align*}
	Y_{r,\kappa}^{(c)}:=d_{r}^{-1}\int_{\Delta_{2}(r)}g(x,y)H_{\kappa}(\xi(x,y))dxdy,
	\end{align*}
	as $r\to\infty$, where $g(x,y)$ is a non-random function such that $g(x,y)\neq0$ when $x=y$, and $d_{r}^{-1}$ is a normalising factor.
	
	We define the corresponding additive functional to $Y_{r,\kappa}^{(c)},\ r>0$, by
	\begin{align}\label{eq6}
	Y_{r,\kappa}^{(d)}&:=d_{r}^{-1}\sum_{i=\lfloor ar\rfloor}^{\lceil br\rceil}\sum_{j=f_{r}^{(l)}(i)}^{f_{r}^{(u)}(i)}g(i,j)H_{\kappa}(\xi(i,j)),
	\end{align}
where
	$f_{r}^{(l)}(i):=\lfloor\displaystyle\inf_{x\in\left[i,i+1\right)}f_{l,r}(x)\rfloor$ and $f_{r}^{(u)}(i):=\lceil\displaystyle\sup_{x\in\left[i,i+1\right)}f_{u,r}(x)\rceil$, see Figure~\ref{fig1:case2b}. It is assumed that
\[
\text{for }ar\in(i,i+1):\ \inf_{x\in\left[i,i+1\right)}f_{l,r}(x) =\lfloor\displaystyle\inf_{x\in\left[ar,i+1\right)}f_{l,r}(x)\rfloor,\ \displaystyle\sup_{x\in\left[i,i+1\right)}f_{u,r}(x) =\lceil\displaystyle\sup_{x\in\left[ar,i+1\right)}f_{u,r}(x)\rceil;
\]
\[
\text{for }br\in(i,i+1):\ \inf_{x\in\left[i,i+1\right)}f_{l,r}(x)=\lfloor\displaystyle\inf_{x\in\left(i,br\right]}f_{l,r}(x)\rfloor, \displaystyle\sup_{x\in\left[i,i+1\right)}f_{u,r}(x) =\lceil\displaystyle\sup_{x\in\left(i,br\right]}f_{u,r}(x)\rceil.
\]
\begin{example}\rm
Figure~\ref{fignew1:case2a} visualises a realisation of the long-range dependent Cauchy field $\xi(x,y)$ over the set $\Delta_{2}$ from Figure~\ref{fig1:case2a}. This field satisfies  Assumptions~\ref{ass1} and \ref{ass1}.  The corresponding normal Q-Q plot of $Y^{(d)}_{r,\kappa}$ in Figure~\ref{fignew1:case2b} was obtained by simulating $\xi(x,y)$ $1000$ times for the large value $r=200$ and $g(x,y)\equiv 1$. It is close to the asymptotic distribution and Figure~\ref{fignew1:case2b} shows its departure from the Gaussian distribution.    
	\begin{figure}
		\centering
		\begin{subfigure}[b]{0.5\textwidth}
			\includegraphics[width=0.85\textwidth,height=6.75cm,trim={0cm 0 0 0},clip]{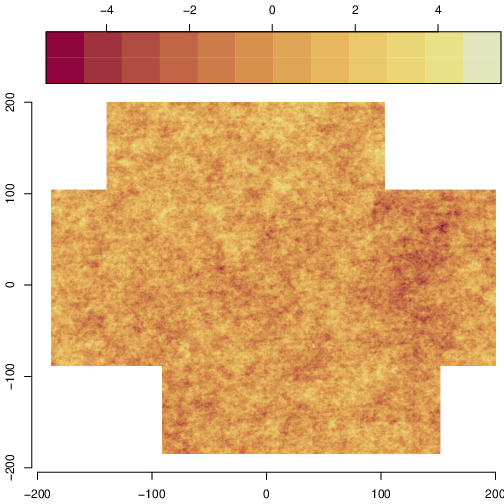}  \vspace{-0.3cm}
			\caption{}
			\label{fignew1:case2a}
		\end{subfigure}
		\hfill  \hspace{-50cm}
		\begin{subfigure}[b]{0.5\textwidth}
			\includegraphics[width=1\textwidth, height=7.7cm,trim={0.5cm 0 1cm 0},clip]{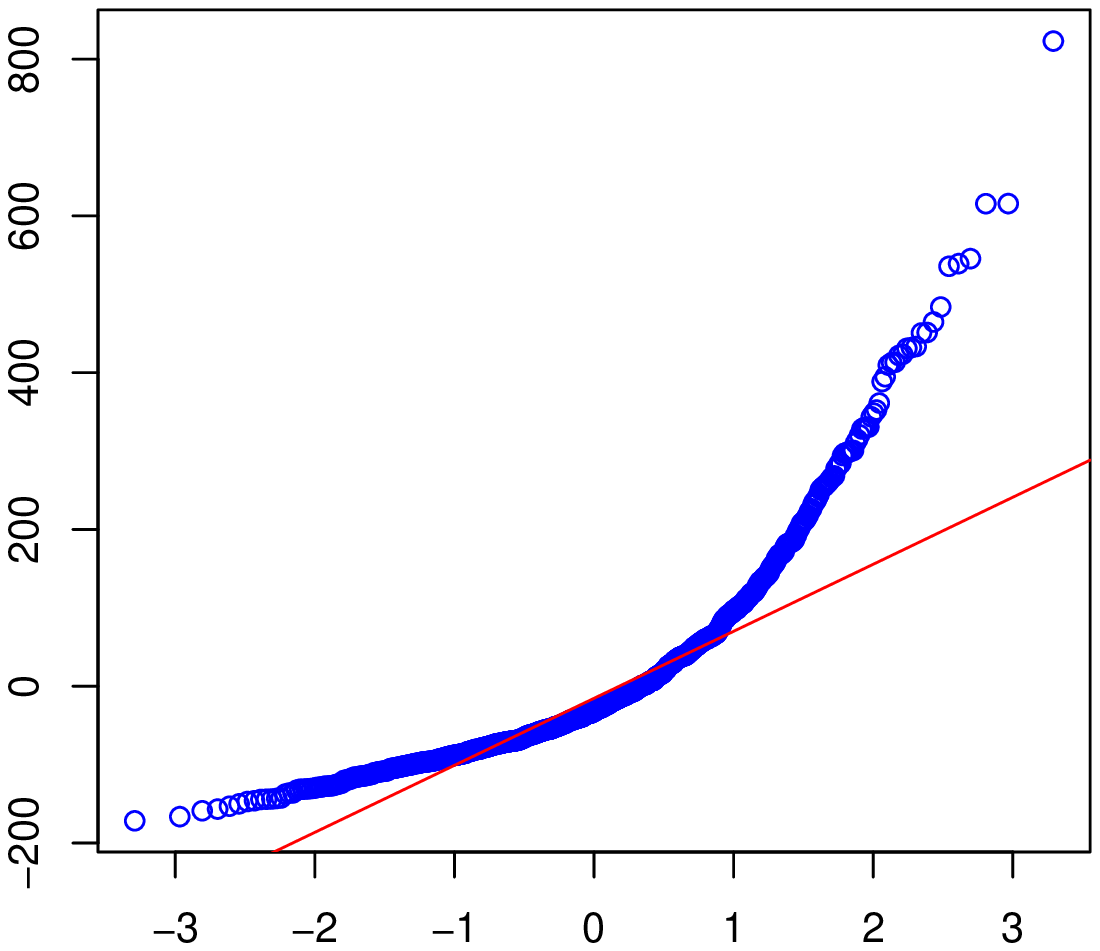} \vspace{-1.7cm}
			\caption{}
			\label{fignew1:case2b}
		\end{subfigure}
		\vspace{-.2cm}
		\caption{(a) A realisation of the Cauchy field $\xi(x,y)$ over $\Delta_{2}$, (b) The normal Q-Q plot of $Y^{(d)}_{200,2}$.}\label{fig4new}
	\end{figure}	
	\end{example}
	\begin{remark}\rm  The functionals $Y_{r,\kappa}^{(c)}$ and $Y_{r,\kappa}^{(d)}$ include various important statistics. For example, the case of $g(x,y)\equiv 1,$ $k=1,$ and $H_1(t)=t$ corresponds to the sample mean estimator.  High order sample moments can be expressed in terms of $Y_{r,\kappa}^{(c)}$ or $Y_{r,\kappa}^{(d)}$ by using the formula
	\[t^\kappa = \kappa!\sum_{m=0}^{\lfloor \kappa/2\rfloor}\frac{1}{2^m m!(\kappa -2m)!}H_{\kappa-2m}(t).\]
	 Another important example, a level excess measure, can be found by  using the Hermite series expansion (\ref{eq2}) for the indicator function $\chi(t > C):$ 
	\[
	\chi(t > C) = \sum\limits_{m = 0}^{\infty}\frac{C_m^{(C)}H_m(t)}{m!},
	\quad  C_m^{(C)} = 
	\begin{cases}
	1 - \Phi(C), & m=0,\\
	\phi(C)H_{m-1}(C), & m\geq 1,
	\end{cases}\]
	where $\Phi(\cdot)$ and $\phi(\cdot)$ are the cdf and pdf for $\mathcal{N}(0,1)$ respectively. The case of $g(x,y)\not\equiv 1$ corresponds to weighted versions of the above statistics.
	\end{remark}

Let us define rectangles $S_{q,r}(i),\ q\in\{l,u\},\ i=\lfloor ar\rfloor,\dots,\lceil br\rceil$, as
\begin{align*}
S_{q,r}(i):=\left[i,i+1\right)\times \bigg[\displaystyle\inf_{x\in\left[i,i+1\right)}f_{q,r}(x),\displaystyle\sup_{x\in\left[i,i+1\right)}f_{q,r}(x)\bigg].
\end{align*}

Then, $Y_{r,\kappa}^{(c)}$ can be rewritten as
\begin{align}\label{eq66}
	Y_{r,\kappa}^{(c)}&=d_{r}^{-1}\int_{A(r)}g(x,y)H_{\kappa}(\xi(x,y))dxdy\notag
	-d_{r}^{-1}\sum_{i=\lfloor ar\rfloor}^{\lceil br\rceil}\int_{S_{l,r}(i)\cap \Delta_{2}^{\mathrm{c}}(r)}g(x,y)H_{\kappa}(\xi(x,y))dydx\notag\\
	&-d_{r}^{-1}\sum_{i=\lfloor ar\rfloor}^{\lceil br\rceil}\int_{S_{u,r}(i)\cap \Delta_{2}^{\mathrm{c}}(r)}g(x,y)H_{\kappa}(\xi(x,y))dydx,
	\end{align}
	where $A(r):=\big\{(x,y)\in\mathbb{R}^{2}:x\in[i,i+1), y\in\big[f_{r}^{(l)}(i),f_{r}^{(u)}(i)\big],\ i=\lfloor ar\rfloor,\dots,\lceil br\rceil\big\}$.
	\begin{assumption}\label{ass4} Let $g(u,v),\ u,v\in\mathbb{R}$, be such that $r^{4-\alpha \kappa}g^{2}\left(r,r \right)L^{\kappa}(r)\to\infty,$ as $r\to\infty,$ and there exists a function $g^{*}(u,v)$ such that for some $\varepsilon>0$ uniformly for $(u,v)\in\Delta_{2\varepsilon}:=\Delta_{2}(1+\varepsilon)$ it holds
	$$\lim_{r\to\infty}\abs{\dfrac{g\left(ru,rv\right)}{g\left(r,r\right)}-g^{*}(u,v)}\rightarrow 0,$$
	where $\alpha\in(0,2/\kappa)$ and
	\begin{align*}
	\int_{\Delta_{2\varepsilon}}\int_{\Delta_{2\varepsilon}}\dfrac{\abs{g^{*}(u_{1},u_{2})g^{*}(v_{1},v_{2})}dv_{1}dv_{2}du_{1}du_{2}}{((u_{1}-v_{1})^{2}+(u_{2}-v_{2})^{2})^{\kappa\alpha/2}}<\infty.
	 \end{align*}	
	\end{assumption}

Notice that, Assumption~\ref{ass4} is more general than Assumption~\ref{ass2weighted}.

	\begin{remark}\label{rem1}\rm It follows from Assumption~{\rm\ref{ass4}} that $g^{*}(u,v)$ is bounded on $\Delta_{2\varepsilon}$.
	\end{remark}
	\begin{remark}\label{rem2}\rm Note, that the conditions on the function $g(\cdot,\cdot)$ in Assumption~{\rm\ref{ass4}} are met by numerous types of functions that are important in solving various statistical problems, in particular, non-linear regression and M estimators. For example, the functions
	\begin{itemize}
		\item $g(u,v)=u^{\mu_{1}}v^{\mu_{2}}\ \text{with}\ g^{*}(u,v)=u^{\mu_{1}}v^{\mu_{2}},$
		\item $g(u,v)=uv\log(\mu_{1}+u)\log(\mu_{2}+v)\ \text{with}\ g^{*}(u,v)=uv$
	\end{itemize}
	(for some appropriate constants $\mu_{1}$ and $\mu_{2}$) can be chosen. The case of $g(u,v)\equiv C>0$ corresponds to the classical equally-weighted functionals and non-central limit theorems.
	\end{remark}
	\begin{remark}\label{rem3}\rm To avoid degenerated cases, the condition $r^{4-\alpha \kappa}g^{2}\left(r,r\right)L^{\kappa}(r)\to\infty,$ as $r\to\infty,$ is essential to guarantee the boundedness of the variance of $d_{r}^{-1}Y_{r,\kappa}^{(c)}$.
	\end{remark}

Now, we proceed to the main result.
	\begin{theorem}\label{theo4} Let $0<\alpha<2/\kappa$. If Assumptions {\rm\ref{ass1}, \ref{ass2}} and {\rm\ref{ass4}} hold, then
	\begin{align}\label{eq7}
		\lim_{\substack{r\to\infty }}\dfrac{\mathbb{E}\left[\int_{\Delta_{2}(r)}g(x,y)H_{\kappa}(\xi(x,y))dydx-\sum_{i=\lfloor ar\rfloor}^{\lceil br\rceil}\sum_{j=f_{r}^{(l)}(i)}^{f_{r}^{(u)}(i)}g(i,j)H_{\kappa}(\xi(i,j))\right]^{2}}{r^{4-\alpha \kappa}L^{\kappa}(r) g^{2}\left(r,r\right)}=0.
		\end{align}
	\end{theorem}
\begin{remark}\label{rem3new}\rm Theorem~\ref{theo4} is also true if $f_{q,1}(\cdot),\ q\in\{l,u\}$, are Lipschitz functions.
\end{remark}
 \subsection{Multidimensional Case}\label{sec5}
This section gives a multidimensional version of Theorem~\ref{theo4} and a generalisation of Theorem~\ref{theo3}. Then, we apply the obtained results to show that additive functionals have the same asymptotic distribution as the corresponding integral functionals. 
 
We use the following notations that enable us to obtain the result of this section analogously to the two-dimensional case. Let $x:=(x_{1},\dots,x_{n})\in\mathbb{R}^{n},\ \textit{1}_{n}:=(1,\dots,1)\in\mathbb{R}^{n}$, and the set $\Delta_{n}\subset\mathbb{R}^{n},\ n\geq3$. 

First, we consider the case of $n=3$. We assume that the set $\Delta_{3}$ can be represented as
\begin{align*}
\Delta_{3}=\{(x_{1},x_{2},x_{3})\in\mathbb{R}^{3}:(x_{1},x_{2})\in\Delta_{2},\ f_{l,1}(x_{1},x_{2})\leq x_{3}\leq f_{u,1}(x_{1},x_{2})\},
\end{align*}
where $f_{l,1}(x_{1},x_{2})<f_{u,1}(x_{1},x_{2})$, $(x_{1},x_{2})\in\Delta_{2}^{\circ}$, and $f_{q,1}(x_{1},x_{2}),\ q\in\set{l,u}$, are smooth functions except of the sets where these functions have finite jumps, i.e. 
\begin{align*}
M_{q}^{\prime}:=\set{(x_{1},x_{2})\in\Delta_{2}:x_{1}=x_{j_{1}}^{q}\ \text{or}\ x_{2}=x_{j_{2}}^{q},\ j_{1}=1,\dots,k_{q}^{(1)},\ j_{2}=1,\dots,k_{q}^{(2)}},
\end{align*} 
where $x_{j_{1}}^{q}$ and $x_{j_{2}}^{q}$, are constants.
Here $k_{q}^{(1)}$ and $k_{q}^{(2)}$ are the number of jumps of $f_{q,1}(\cdot,\cdot)$. Thus, $M_{q}^{\prime}$ consists of a finite number of two-dimensional line segments. 

Note that, by the  homothety of $\Delta_{3}(r),\ r>0$, the set $\Delta_{3}(r)$ can be represented as $\Delta_{3}(r)=\{(x_{1},x_{2},x_{3})\in\mathbb{R}^{3}:(x_{1},x_{2})\in\Delta_{2}(r),\ f_{l,r}(x_{1},x_{2})\leq x_{3}\leq f_{u,r}(x_{1},x_{2})\}$, where $f_{q,r}(x_{1},x_{2})=rf_{q,1}(x_{1}/r,x_{2}/r)$, $q\in\set{l,u}$.

For a real-valued homogeneous isotropic Gaussian random field $\xi(x),\ x\in\mathbb{R}^{3}$, let
\begin{align*}
Z_{r,\kappa}^{(c)}:=d_{r}^{-1}\int_{\Delta_{3}(r)}g(x)H_{\kappa}(\xi(x))dx,
\end{align*}
where $g(x),\ x\in\mathbb{R}^{3}$, is a non-random scalar function such that $g(x_{1}\textit{1}_{3})\neq0$, and $d_{r}^{-1}$ is a normalising factor.

We define the corresponding additive functional to $Z_{r,\kappa}^{(c)},\ r>0$, by
\begin{align*}
Z_{r,\kappa}^{(d)}&:=d_{r}^{-1}\sum_{(i_{1},i_{2},i_{3})\in Q_{3}(\Delta_{3}(r))}g(i_{1},i_{2},i_{3})H_{\kappa}(\xi(i_{1},i_{2},i_{3})),
\end{align*}
where
\begin{align*}
Q_{3}(\Delta_{3}(r)):=\bigg\{&(i_{1},i_{2},i_{3})\in\mathbb{Z}^{3}:\ i_{1}\in\set{\lfloor ar\rfloor,\dots,\lceil br\rceil},\ i_{2}\in\set{f_{r}^{(l)}(i_{1}),\dots,f_{r}^{(u)}(i_{1})},\\ &i_{3}\in\set{f_{r}^{(l)}(i_{1},i_{2}),\dots,f_{r}^{(u)}(i_{1},i_{2})}\bigg\},
\end{align*}
\begin{align*}
f_{r}^{(l)}(i_{1},i_{2}):=\bigg\lfloor\inf_{\substack{(x_{1},x_{2})\in P_{2}(r,i_{1},i_{2})}}f_{l,r}(x_{1},x_{2})\bigg\rfloor\ \text{and}\  f_{r}^{(u)}(i_{1},i_{2}):=\bigg\lfloor\sup_{\substack{(x_{1},x_{2})\in P_{2}(r,i_{1},i_{2})}}f_{u,r}(x_{1},x_{2})\bigg\rfloor,	
\end{align*}
and
\begin{align*}
	P_{2}(r,i_{1},i_{2}):=\bigg\{&(x_{1},x_{2})\in\left[i_{1},i_{1}+1\right)\times\left[i_{2},i_{2}+1\right]\cap\Delta_{2}(r),\ i_{1}\in\set{\lfloor ar\rfloor,\dots,\lceil br\rceil},\\ &i_{2}\in[f^{(l)}_{r}(i_{1}),f^{(u)}_{r}(i_{1})]\bigg\}.
\end{align*}
\begin{remark}\rm
	The intersection with $\Delta_{2}(r)$ is required to correctly define the infimum and supremum for the cases when some points in $\left[i_{1},i_{1}+1\right)\times\left[i_{2},i_{2}+1\right]$ are outside of $\Delta_{2}(r)$, which may happen for the boundary region. 
\end{remark}
For each $q\in\{l,u\}$, define three-dimensional parallelepipeds $S_{q,r}(i_{1},i_{2})$ as
\[
S_{q,r}(i_{1},i_{2}):=\left[i_{1},i_{1}+1\right)\times\left[i_{2},i_{2}+1\right]\times\left[\inf_{\substack{(x_{1},x_{2})\in P_{2}(r,i_{1},i_{2})}}f_{q,r}(x_{1},x_{2}),\sup_{\substack{(x_{1},x_{2})\in P_{2}(r,i_{1},i_{2})}}f_{q,r}(x_{1},x_{2})\right],
\]
where $i_{1}=\lfloor ar\rfloor,\dots,\lceil br\rceil,\ i_{2}=f^{(l)}_{r}(i_{1}),\dots,f^{(u)}_{r}(i_{1})$.

By induction one can extend this construction to an arbitrary dimension $n$ as follows.

The sets $\Delta_{n}$ and $\Delta_{n}(r),\ n\geq4$, can be defined as
\begin{align}\label{eqA1}
\Delta_{n}=\{(x_{1},\dots,x_{n-1})\in\Delta_{n-1},\ f_{l,1}(x_{1},\dots,x_{n-1})\leq x_{n}\leq f_{u,1}(x_{1},\dots,x_{n-1})\},
\end{align}
and
$\Delta_{n}(r)=\{(x_{1},\cdots,x_{n-1})\in\Delta_{n-1}(r), f_{l,r}(x_{1},\dots,x_{n-1})\leq x_{n}\leq f_{u,r}(x_{1},\dots,x_{n-1})\},$ 
where
$f_{q,r}(x_{1},\dots,x_{n-1})=rf_{q,1}(\frac{x_{1}}{r},\dots,\frac{x_{n-1}}{r})$, such that $f_{l,1}(x_{1},\dots,x_{n-1})<f_{u,1}(x_{1},\dots,x_{n-1})$ if $(x_{1},\dots,x_{n-1})\in\Delta_{n-1}^{\circ}$, and $f_{q,1}(x_{1},\dots,x_{n-1}),$ $\ q\in\set{l,u}$, are smooth functions except of the sets where these functions have finite jumps, i.e. 
\begin{align}\label{eqA2}
M_{q}^{\prime\prime}:=\set{(x_{1},\dots,x_{n-1})\in\Delta_{n-1}:x_{1}=x_{j_{1}}^{q}\ \text{or}\ x_{2}=x_{j_{2}}^{q}\ \text{or}\cdots\ \text{or}\ x_{n-1}=x_{j_{n-1}}^{q}},
\end{align} 
where $j_{m}=1,\dots,k_{q}^{(m)},\ m=1,\dots,n-1$, and $x_{j_{m}}^{q},\ m=1,\dots,n-1$, are constants.
Here $k_{q}^{(m)},\ m=1,\dots,n-1$, is the number of jumps of $f_{q,1}(\cdot,\dots,\cdot)$ over dimension $m$. Thus, $M_{q}^{\prime\prime}$ consists of a finite number of $(n-1)$-dimensional hyperplanes sections in $\Delta_{n-1}$.

For a real-valued homogeneous isotropic Gaussian random field $\xi(x),\ x\in\mathbb{R}^{n}$, let
\begin{align*}
Z_{r,\kappa}^{(c)}:=d_{r}^{-1}\int_{\Delta_{n}(r)}g(x)H_{\kappa}(\xi(x))dx,
\end{align*}
where $g(x),\ x\in\mathbb{R}^{n}$, is a non-random function such that $g(x_{1}\textit{1}_{n})\neq0$, and $d_{r}^{-1}$ is a normalising factor.

We define the corresponding additive functional to $Z_{r,\kappa}^{(c)},\ r>0$, by
\begin{align*}
Z_{r,\kappa}^{(d)}&:=d_{r}^{-1}\sum_{\textit{\textbf{i}}\in Q_{n}(\Delta_{n}(r))}g(\textit{\textbf{i}})H_{\kappa}(\xi(\textit{\textbf{i}})),
\end{align*}
where $\textit{\textbf{i}}=(i_{1},\dots,i_{n})\in\mathbb{Z}^{n}$,
\begin{align*}
Q_{n}(\Delta_{n}(r)):=\bigg\{\textit{\textbf{i}}&:(i_{1},\dots,i_{n-1})\in Q_{n-1}(\Delta_{n-1}(r)),\\ &i_{n}\in\set{f_{r}^{(l)}(i_{1},\dots,i_{n-1}),\dots,f_{r}^{(u)}(i_{1},\dots,i_{n-1})}\bigg\},
\end{align*}
\[
f_{r}^{(l)}(i_{1},\dots,i_{n-1}):=\bigg\lfloor\inf_{\substack{(x_{1},\dots,x_{n-1})\in P_{n-1}(r,i_{1},\dots,i_{n-1})}}f_{l,r}(x_{1},\dots,x_{n-1})\bigg\rfloor,
\]
\[
f_{r}^{(u)}(i_{1},\dots,i_{n-1}):=\bigg\lfloor\sup_{\substack{(x_{1},\dots,x_{n-1})\in P_{n-1}(r,i_{1},\dots,i_{n-1})}}f_{u,r}(x_{1},\dots,x_{n-1})\bigg\rfloor,
\]
and
\begin{align*}
&P_{n-1}(r,i_{1},\dots,i_{n-1}):=\bigg\{ (x_{1},\dots,x_{n-1})\in \left[i_{1},i_{1}+1\right)\times\left[i_{2},i_{2}+1\right]\times\cdots\times\left[i_{n-1},i_{n-1}+1\right]\cap\\
&\Delta_{n-1}(r), i_{1}=\lfloor ar\rfloor,\dots,\lceil br\rceil,\ i_{j}\in[f^{(l)}_{r}(i_{1},\dots,i_{j-1}),f^{(u)}_{r}(i_{1},\dots,i_{j-1})],\ j=2,\dots,n-1\bigg\}.
\end{align*}
\begin{assumption}\label{ass5} Let $g(u),\ u\in\mathbb{R}^{n}$, be such that $r^{2n-\alpha \kappa}g^{2}\left(r\textit{1}_{n} \right)L^{\kappa}(r)\rightarrow\infty,$ as $r\to\infty,$ and there exists a function $g^{*}(u)$ such that for some $\varepsilon>0$ uniformly for $u\in\Delta_{n\varepsilon}:=\Delta_{n}(1+\varepsilon)$ it holds
	$$\lim_{r\to\infty}\abs{\dfrac{g\left(ru\right)}{g\left(r\textit{1}_{n}\right)}-g^{*}(u)}\rightarrow 0$$
	where $\alpha\in(0,n/\kappa)$ and
	\begin{align*}
	\int_{\Delta_{n\varepsilon}}\int_{\Delta_{n\varepsilon}}\dfrac{\abs{g^{*}(u)g^{*}(v)}dudv}{\Vert u-v\Vert^{\kappa\alpha}}<\infty.
	\end{align*}	
\end{assumption}

Following the steps analogous to the proof in Section~\ref{secproofs} and replacing intervals by multidimensional parallelepipeds we obtain a multidimensional version of Theorem~\ref{theo4}.

\begin{theorem}\label{theo5} Let $\xi(x),\ x\in\mathbb{R}^{n}$, and $\Delta_{n}$ satisfies assumptions {\rm(\ref{eqA1})} and {\rm(\ref{eqA2})}. If Assumptions~{\rm\ref{ass1}, \ref{ass2}} and {\rm\ref{ass5}} hold, $\alpha\in\left(0,n/\kappa\right)$, then
	\begin{equation*}\label{eq21}
	\lim_{\substack{r\to\infty }}\dfrac{\mathbb{E}\left[\int_{\Delta_{n}(r)}g(x)H_{\kappa}(\xi(x))dx-\sum_{\textit{\textbf{i}}\in Q_{n}(\Delta_{n}(r))}g(\textit{\textbf{i}})H_{\kappa}(\xi(\textit{\textbf{i}}))\right]^{2}}{r^{2n-\alpha\kappa}L^{\kappa}(r) g^{2}\left(r\textit{1}_{n}\right)}=0.
	\end{equation*}
\end{theorem}
 
Next, we give a generalisation of Theorem~\ref{theo3}.
\begin{assumption}\label{ass6}
Let $g^{*}(x),\ x\in\mathbb{R}^{n}$, be a function such that for $\alpha\in\left(0,n/\kappa\right)$ it holds 
\begin{align*}
\int_{\mathbb{R}^{n\kappa}}\prod_{j=1}^{\kappa}\Vert \lambda_{j}\Vert^{\alpha-n}\abs{K_{\Delta_{n}}\left(\lambda;g^{*}\right)}^{2}\prod_{j=1}^{\kappa}d\lambda_{j}<\infty
\end{align*}
and
\begin{align*}
\lim_{r\to\infty}\int_{\mathbb{R}^{n\kappa}}\bigg|\int_{\Delta_{n}}e^{i\langle \lambda_{1}+\dots+\lambda_{\kappa},x \rangle}
	\bigg(\frac{g(rx)}{g(r\textit{1}_{n})}\prod_{j=1}^{\kappa}\sqrt{\frac{L(r/\Vert\lambda_{j}\Vert)}{L(r)}}-g^{*}(x)\bigg)dx\bigg|^{2}\prod_{j=1}^{\kappa}\Vert \lambda_{j}\Vert^{\alpha-n}\prod_{j=1}^{\kappa}d\lambda_{j}=0.
\end{align*}
\end{assumption}

\begin{theorem}\label{theonew7} Let $\Delta_{n}(r)$ satisfy~{\rm(\ref{eqA1})} and~{\rm(\ref{eqA2})}. If Assumptions~{\rm\ref{ass1}}, {\rm\ref{ass2}}, {\rm\ref{ass5}} and {\rm\ref{ass6}} hold, $\alpha\in\left(0,\frac{n}{\kappa}\right)$, then for $r\to\infty$ 
	\begin{equation*}\label{tt}
		Z^{(c)}_{r,\kappa}:=\dfrac{1}{r^{n-\alpha \kappa/2} g\left(r\textit{1}_{n}\right)L^{\kappa/2}(r)c_{1}^{\kappa/2}(n,\alpha)}\int_{\Delta_{n}(r)}g(x)H_{\kappa}(\xi(x))dx
	\end{equation*}
	converge weakly to the random variable
	\begin{align*}
	Z_{\kappa}:=\int_{\mathbb{R}^{n\kappa}}^{\prime}K_{\Delta_{n}}\left(\lambda_{1}+\cdots+\lambda_{\kappa};g^{*}\right)\dfrac{\prod_{j=1}^{\kappa} W(d\lambda_{j})}{\prod_{j=1}^{\kappa}\Vert\lambda_{j}\Vert^{(n-\alpha)/2}}.
	\end{align*}
	\end{theorem}

Now we apply the result of Theorem~\ref{theo5} to Theorem~\ref{theonew7} to obtain an analogous result in the discrete case.
\begin{theorem}\label{theo7} Let $\Delta_{n}(r)$ satisfy~{\rm(\ref{eqA1})} and {\rm(\ref{eqA2})}. If Assumptions~{\rm\ref{ass1}},~{\rm\ref{ass2}}, {\rm\ref{ass5}} and {\rm\ref{ass6}} hold, $\alpha\in\left(0,\frac{n}{\kappa}\right)$, then for $r\to\infty$
	\begin{equation*}
	Z^{(d)}_{r,\kappa}:=\dfrac{1}{r^{n-\alpha \kappa/2}L^{\kappa/2}(r) g\left(r\textit{1}_{n}\right)c_{1}^{\kappa/2}(n,\alpha)}\sum_{\textit{\textbf{i}}\in Q_{n}(\Delta_{n}(r))}g(\textit{\textbf{i}})H_{\kappa}(\xi(\textit{\textbf{i}}))
	\end{equation*}
	converge weakly to the random variable $Z_{\kappa}.$
	
\begin{remark}\rm	 If  $\kappa= 1$ the limit $Z_{\kappa}$ is Gaussian. For 
	$\kappa > 1$  the random variables $Z_{\kappa}$ have non-Gaussian distribution. The most studied case is the Rosenblatt distribution that corresponds $\kappa = 2$ and a rectangular $\Delta_{n},$ see~\cite{taqqu2013}.
\end{remark}
	
\end{theorem}
\section{Proofs of Results from Sections~\ref{sec4} and~\ref{sec5} }\label{secproofs}
In this section we give proofs of results in Sections~\ref{sec4} and~\ref{sec5}.
\begin{flushleft}
	\textit{Proof of Theorem~\rm{\ref{theo4}}}. Using (\ref{eq6}) and (\ref{eq66}), one can estimate the numerator in (\ref{eq7}) as
\end{flushleft}
\begin{align}
&\mathbb{E}\bigg[\int_{A(r)}g(x,y)H_{\kappa}(\xi(x,y))dydx-\sum_{i=\lfloor ar\rfloor}^{\lceil br\rceil}\sum_{j=f_{r}^{(l)}(i)}^{f_{r}^{(u)}(i)}g(i,j)H_{\kappa}(\xi(i,j))\notag\\
&-\sum_{i=\lfloor ar\rfloor}^{\lceil br\rceil}\int_{S_{l,r}(i)\cap \Delta_{2}^{\mathrm{c}}(r)}g(x,y)H_{\kappa}(\xi(x,y))dydx-\sum_{i=\lfloor ar\rfloor}^{\lceil br\rceil}\int_{S_{u,r}(i)\cap \Delta_{2}^{\mathrm{c}}(r)}g(x,y)H_{\kappa}(\xi(x,y))dydx\bigg]^{2}\notag\\
&\leq 2\mathbb{E}\bigg[\int_{A(r)}g(x,y)H_{\kappa}(\xi(x,y))dydx-\sum_{i=\lfloor ar\rfloor}^{\lceil br\rceil}\sum_{j=f_{r}^{(l)}(i)}^{f_{r}^{(u)}(i)}g(i,j)H_{\kappa}(\xi(i,j))\bigg]^{2}\notag\\
&+2\mathbb{E}\bigg[\sum_{i=\lfloor ar\rfloor}^{\lceil br\rceil}\bigg(\int_{S_{l,r}(i)\cap \Delta_{2}^{\mathrm{c}}(r)}g(x,y)H_{\kappa}(\xi(x,y))dydx+\int_{S_{u,r}(i)\cap \Delta_{2}^{\mathrm{c}}(r)}g(x,y)H_{\kappa}(\xi(x,y))dydx\bigg)\bigg]^{2}\notag\\
&\leq 2\mathbb{E}\bigg[\int_{A(r)}g(x,y)H_{\kappa}(\xi(x,y))dydx-\sum_{i=\lfloor ar\rfloor}^{\lceil br\rceil}\sum_{j=f_{r}^{(l)}(i)}^{f_{r}^{(u)}(i)}g(i,j)H_{\kappa}(\xi(i,j))\bigg]^{2}\notag\\
&+4\mathbb{E}\bigg[\sum_{i=\lfloor ar\rfloor}^{\lceil br\rceil}\int_{S_{l,r}(i)\cap \Delta_{2}^{\mathrm{c}}(r)}g(x,y)H_{\kappa}(\xi(x,y))dydx\bigg]^{2}\notag\\
&+4\mathbb{E}\bigg[\sum_{i=\lfloor ar\rfloor}^{\lceil br\rceil}\int_{S_{u,r}(i)\cap \Delta_{2}^{\mathrm{c}}(r)}g(x,y)H_{\kappa}(\xi(x,y))dydx\bigg]^{2}=: 2J_{1}+4J_{2}+4J_{3}.\label{eq9}
\end{align}
We will consider each term in (\ref{eq9}) separately. By~(\ref{eq1}) we get

\begin{align*}
&\dfrac{J_{2}}{r^{4-\kappa\alpha}g^{2}(r,r)L^{\kappa}(r)}=\dfrac{\mathbb{E}\bigg[\sum_{i=\lfloor ar\rfloor}^{\lceil br\rceil}\int_{S_{l,r}(i)\cap \Delta_{2}^{\mathrm{c}}(r)}g(x,y)H_{\kappa}(\xi(x,y))dydx\bigg]^{2}}{r^{4-\kappa\alpha}g^{2}(r,r)L^{\kappa}(r)}\\
&=\kappa!\sum_{i,j=\lfloor ar\rfloor}^{\lceil br\rceil}\int_{S_{l,r}(i)\cap \Delta_{2}^{\mathrm{c}}(r)}\int_{S_{l,r}(j)\cap \Delta_{2}^{\mathrm{c}}(r)}\dfrac{g(x,y)g(x',y')B^{\kappa}(\Vert(x-x',y-y')\Vert)dydy'dxdx'}{r^{4-\kappa\alpha}g^{2}(r,r)L^{\kappa}(r)}\\
&\leq \kappa!\sum_{i,j=\lfloor ar\rfloor}^{\lceil br\rceil}\int_{S_{l,r}(i)}\int_{S_{l,r}(j)}\dfrac{\vert g(x,y)g(x',y')B^{\kappa}(\Vert(x-x',y-y')\Vert)\vert dydy'dxdx'}{r^{4-\kappa\alpha}g^{2}(r,r)L^{\kappa}(r)}.
\end{align*}
As by Assumption~\ref{ass1} $\vert B(\cdot)\vert\leq 1$, we have
\begin{align*}	
\dfrac{J_{2}}{r^{4-\kappa\alpha}g^{2}(r,r)L^{\kappa}(r)}\leq\dfrac{\kappa!}{r^{4-\kappa\alpha}L^{\kappa}(r)}\sum_{i,j=\lfloor ar\rfloor}^{\lceil br\rceil}\int_{S_{l,r}(i)}\int_{S_{l,r}(j)}\dfrac{\vert g(x,y)g(x',y')\vert dydy'dxdx'}{g^{2}(r,r)}.
\end{align*}
Using the following transformation
\begin{align}\label{eq15}
ru_{1}=x,ru_{2}=y, rv_{1}=x',\ \textrm{and}\ rv_{2}=y',
\end{align}
we get
\begin{align}\label{eq10}
\dfrac{J_{2}}{r^{4-\kappa\alpha}g^{2}(r,r)L^{\kappa}(r)}\leq\dfrac{\kappa!}{r^{-\kappa\alpha}L^{\kappa}(r)}\sum_{i,j=\lfloor ar\rfloor}^{\lceil br\rceil}\int_{S_{l,r}^{*}(i)}\int_{S_{l,r}^{*}(j)}\dfrac{\vert g(ru_{1},ru_{2})g(rv_{1},rv_{2})\vert du_{1}du_{2}dv_{1}dv_{2}}{g^{2}(r,r)},
\end{align}
where
\begin{align*}
S_{l,r}^{*}(i)=\set{(u_{1},u_{2}):u_{1}\in\left[\frac{i}{r},\frac{i+1}{r}\right),\ u_{2}\in\bigg[\displaystyle\inf_{u_{1}\in\left[\frac{i}{r},\frac{i+1}{r}\right)}f_{l,1}(u_{1}),\displaystyle\sup_{u_{1}\in\left[\frac{i}{r},\frac{i+1}{r}\right)}f_{l,1}(u_{1})\bigg)}.
\end{align*}

Note, that the rectangles $\{S_{l,r}^{*}(i),i\in \left[\lfloor ar\rfloor,\lceil br\rceil\right]\}$ have the same width but different lengths. Denote by $i'\in\{\lfloor ar\rfloor,\dots,\lceil br\rceil\}$ a such index that the rectangle $S_{l,r}^{*}(i')$ has the largest length. Let $\varepsilon_{*}:=\min_{\substack{i,j=1,\dots,k_{l}\\ i\neq j}}\vert x^{l}_{i}-x^{l}_{j}\vert$ and for each point of discontinuity $x^{l}_{j}\in\ M_{l}$ the $\varepsilon_{*}/r$ neighbourhood of $x^{l}_{j}$ be defined by $N_{\varepsilon_{*}}(x^{l}_{j},r):=\left(x^{l}_{j}-\frac{\varepsilon_{*}}{r},x^{l}_{j}+\frac{\varepsilon_{*}}{r}\right)$. Let $N_{\varepsilon_{*}}^{\prime}(x^{l}_{j},r)=[x^{l}_{j}-\frac{\varepsilon_{*}}{r},x^{l}_{j}+\frac{\varepsilon_{*}}{r}]\setminus\{x^{l}_{j}\}$ and
\begin{align*}
SN_{1}(x^{l}_{j},r):&=N_{\epsilon_{*}}(x^{l}_{j},r)\times\bigg(
\inf_{u_{1}\in N_{\epsilon_{*}}^{\prime}(x^{l}_{j},r)}f_{l,1}(u_{1}) ,\displaystyle\sup_{u_{1}\in N_{\epsilon_{*}}^{\prime}(x^{l}_{j},r)}f_{l,1}(u_{1})\bigg).
\end{align*}
We also define
\begin{align*}
T_{k_{l}}:&=\set{i\in\set{\lfloor ar\rfloor,\dots,\lceil br\rceil}:x^{l}_{j}\not\in\left[\frac{i}{r},\frac{i+1}{r}\right)\ \text{for all}\  j=1,\dots,k_{l}},
\end{align*}
and for $x^{l}_{j}\in[\frac{i^{\prime}}{r},\frac{i^{\prime}+1}{r})$
\begin{align*}
SN_{2}(x^{l}_{j},r):&=\left[\frac{i^{\prime}}{r},x^{l}_{j}-\frac{\varepsilon_{*}}{r}\right)\times\bigg[\displaystyle\inf_{u_{1}\in[\frac{i^{\prime}}{r},x^{l}_{j}-\frac{\varepsilon_{*}}{r})}f_{l,1}(u_{1}),\displaystyle\sup_{u_{1}\in[\frac{i^{\prime}}{r},x^{l}_{j}-\frac{\varepsilon_{*}}{r})}f_{l,1}(u_{1}])\bigg],\\
SN_{3}(x^{l}_{j},r):&=\left(x^{l}_{j}+\frac{\varepsilon_{*}}{r},\frac{i^{\prime}+1}{r}\right)\times\bigg[\displaystyle\inf_{u_{1}\in(x^{l}_{j}+\frac{\varepsilon_{*}}{r},\frac{i^{\prime}+1}{r})}f_{l,1}(u_{1}),\displaystyle\sup_{u_{1}\in(x^{l}_{j}+\frac{\varepsilon_{*}}{r},\frac{i^{\prime}+1}{r})}f_{l,1}(u_{1})\bigg].
\end{align*}

Note, that for each $j=1,\dots,k_{l}$, $f_{l,1}(\cdot)$ is a bounded function on $N_{\epsilon_{*}}(x^{l}_{j},r)$. As the number of jumps is finite then there is a constant $C>0$, such that for all $j=1,\dots,k_{l}$ it holds
\begin{align*}
\big\vert SN_{1}(x^{l}_{j},r)\big\vert&\leq 2\dfrac{\varepsilon_{*}}{r}\bigg\vert\displaystyle\sup_{u_{1}\in N_{\varepsilon_{*}}(x^{l}_{j},r)}f_{l,1}(u_{1})-\inf_{u_{1}\in N_{\varepsilon_{*}}(x^{l}_{j},r)}f_{l,1}(u_{1})\bigg\vert\leq C\frac{\varepsilon_{*}}{r}.
\end{align*}

The smoothness of the function $f_{l,1}(\cdot)$ in $[\frac{i^{\prime}}{r},x^{l}_{j}-\frac{\varepsilon_{*}}{r})$ and $(x^{l}_{j}+\frac{\varepsilon_{*}}{r},\frac{i^{\prime}+1}{r}]$ gives
\begin{align*}
\big\vert SN_{e}(x^{l}_{j},r)\big\vert\leq \frac{C}{r},\quad e=2,3.
\end{align*}

Now, using the above results for sufficient large $r$, one can estimate~(\ref{eq10}) as
\begin{align*}
\dfrac{J_{2}}{r^{4-\kappa\alpha}g^{2}(r,r)L^{\kappa}(r)}&\leq\dfrac{\kappa!}{r^{-\kappa\alpha}L^{\kappa}(r)}\bigg(\sup_{i\in\left[\lfloor ar\rfloor,\lceil br\rceil\right]}\sup_{\substack{(u_{1},u_{2})\in S^{*}_{l,r}(i)\\}}\bigg\vert\dfrac{g(ru_{1},ru_{2})}{g(r,r)}\bigg\vert\bigg)^{2}\quad\quad\quad\quad
\end{align*}
\begin{align*}
&\times\bigg(\left(\lceil br\rceil-\lfloor ar\rfloor+1\right)\sup_{\substack{i\in T_{k_{l}}}}\vert S_{l,r}^{*}(i)\vert+k_{l}\sum_{e=1}^{3}\sup_{\substack{j=1,\cdots,k_{l}}}\big\vert SN_{e}(x^{l}_{j},r)\big\vert\bigg)^{2}\\
&\leq\dfrac{\kappa!}{r^{-\kappa\alpha}L^{\kappa}(r)}\bigg(\sup_{i\in\left[\lfloor ar\rfloor,\lceil br\rceil\right]}\sup_{\substack{(u_{1},u_{2})\in S^{*}_{l,r}(i)\\}}\bigg\vert\dfrac{g(ru_{1},ru_{2})}{g(r,r)}\bigg\vert\bigg)^{2}\notag\\
&\times\left(\frac{\left(\lceil br\rceil-\lfloor ar\rfloor+1\right)}{r}\bigg\vert \displaystyle\sup_{u_{1}\in\left[\frac{i^{\prime}}{r},\frac{i^{\prime}+1}{r}\right)}f_{l,1}(u_{1})-\displaystyle\inf_{u_{1}\in\left[\frac{i^{\prime}}{r},\frac{i^{\prime}+1}{r}\right)}f_{l,1}(u_{1})\bigg\vert+\dfrac{k_{l}C}{r}\right)^{2}
\end{align*}
\begin{align}\label{eq11}
&\leq \dfrac{C}{r^{2-\kappa\alpha}L^{\kappa}(r)}\bigg(\sup_{i\in\left[\lfloor ar\rfloor,\lceil br\rceil\right]}\sup_{\substack{(u_{1},u_{2})\in S^{*}_{l,r}(i)\\}}\bigg\vert\dfrac{g(ru_{1},ru_{2})}{g(r,r)}\bigg\vert\bigg)^{2}\notag\\
&\times\left(\left(\lceil br\rceil-\lfloor ar\rfloor+1\right)\bigg\vert \displaystyle\sup_{u_{1}\in\left[\frac{i^{\prime}}{r},\frac{i^{\prime}+1}{r}\right)}f_{l,1}(u_{1})-\displaystyle\inf_{u_{1}\in\left[\frac{i^{\prime}}{r},\frac{i^{\prime}+1}{r}\right)}f_{l,1}(u_{1})\bigg\vert+C\right)^{2}.
\end{align}

Note, that as the function $f_{l,1}(\cdot)\in\mathcal{C}^{1}$ on $T_{k_{l}}$, then by the mean-value theorem there exists  $u_{0}\in\left[\frac{i}{r},\frac{i+1}{r}\right),\ i\in T_{k_{l}}$, such that
\begin{align*}
\bigg\vert\displaystyle\sup_{u_{1}\in\left[\frac{i}{r},\frac{i+1}{r}\right)}f_{l,1}(u_{1})-\displaystyle\inf_{u_{1}\in\left[\frac{i}{r},\frac{i+1}{r}\right)}f_{l,1}(u_{1})\bigg\vert\leq \frac{1}{r}\max_{i\in T_{k_{l}}}\sup_{u_{0}\in\left[\frac{i}{r},\frac{i+1}{r}\right)}\big\vert f_{l,1}^{\prime}(u_{0})\big\vert=\frac{\tilde{C}}{r}.
\end{align*}

Then,~(\ref{eq11}) can be estimated as
\begin{align}\label{eq11new}
\dfrac{J_{2}}{r^{4-\kappa\alpha}g^{2}(r,r)L^{\kappa}(r)}&\leq \dfrac{\tilde{\tilde{C}}}{r^{2-\kappa\alpha}L^{\kappa}(r)}\bigg(\sup_{i\in\left[\lfloor ar\rfloor,\lceil br\rceil\right]}\sup_{\substack{(u_{1},u_{2})\in S^{*}_{l,r}(i)\\}}\bigg\vert\dfrac{g(ru_{1},ru_{2})}{g(r,r)}\bigg\vert\bigg)^{2}\notag\\
&\times\left(\tilde{C}\left(\frac{\lceil br\rceil-\lfloor ar\rfloor+1}{r}\right)+C\right)^{2}.
\end{align} 

As for sufficiently large $r$ the rectangles $S^{*}_{l,r}(i)\subset\Delta_{2\varepsilon}$ for all $i=\lfloor ar\rfloor,\dots,\lceil br\rceil$, then by Assumption~\ref{ass4} and Remark~\ref{rem1} for an arbitrary $\varepsilon_{1}>0$ there exists $r_{0}$ such that for $r>r_{0}$
\begin{align*}
\max_{i=\lfloor ar\rfloor,\dots,\lceil br\rceil}\sup_{\substack{(u_{1},u_{2})\in S^{*}_{l,r}(i)\\}}\bigg\vert\dfrac{g(ru_{1},ru_{2})}{g(r,r)}\bigg\vert\leq\sup_{\substack{(u_{1},u_{2})\in \Delta_{2\varepsilon}\\}}\vert
g^{*}(u_{1},u_{2})\vert+\varepsilon_{1}<\infty.
\end{align*}

As $a=\lim\limits_{r\to\infty}\frac{\lfloor ar\rfloor}{r}$ and $b=\lim\limits_{r\to\infty}\frac{\lceil br\rceil}{r}$ it holds that $\frac{\lceil br\rceil-\lfloor ar\rfloor+1}{r}\to b-a<\infty$, when~$r\to\infty$.

Finally, as $0<\alpha<2/\kappa$ the term $\dfrac{1}{r^{2-\kappa\alpha}L^{\kappa}(r)}\to 0$, when  $r\to\infty$. Hence, the upper bound in (\ref{eq11new}) approaches $0$ when $r\to\infty$.

Similarly, one obtains that $\dfrac{J_{3}}{r^{4-\kappa\alpha}g^{2}(r,r)L^{\kappa}(r)}\to 0, \quad {\rm as}\ r\to 0.$

Therefore, it is enough to investigate the behaviour of
\begin{align}\label{eq12}
J_{1}=\mathbb{E}\bigg[\int_{A(r)}g(x,y)H_{\kappa}(\xi(x,y))dydx-\sum_{i=\lfloor ar\rfloor}^{\lceil br\rceil}\sum_{j=f_{r}^{(l)}(i)}^{f_{r}^{(u)}(i)}g(i,j)H_{\kappa}(\xi(i,j))\bigg]^{2}.
\end{align}

One can estimate (\ref{eq12}) as
\begin{align}\label{eq13}
J_{1}&=\mathbb{E}\bigg(\sum_{i=\lfloor ar\rfloor}^{\lceil br\rceil}\sum_{j=f_{r}^{(l)}(i)}^{f_{r}^{(u)}(i)}\int_{[0,1]^{2}}(g(x+i,y+j)H_{\kappa}(\xi(x+i,y+j))-g(i,j)H_{\kappa}(\xi(i,j)))dydx\bigg)^{2}\notag\\
&=:\sum_{k=1}^{3}D_{r}^{(k)},
\end{align}
where
\begin{align*}    
D_{r}^{(1)}&=\mathbb{E}\sum_{i,i'=\lfloor ar\rfloor}^{\lceil br\rceil}\sum_{j=f_{r}^{(l)}(i)}^{f_{r}^{(u)}(i)}\sum_{j'=f_{r}^{(l)}(i')}^{f_{r}^{(u)}(i')}\int_{\left[0,1\right]^{4}}g(x+i,y+j)g(x'+i',y'+j')H_{\kappa}(\xi(x+i,y+j))\\
&\times H_{\kappa}(\xi(x'+i',y'+j'))dydy'dxdx',\\
D_{r}^{(2)}&=-2\mathbb{E}\sum_{i,i'=\lfloor ar\rfloor}^{\lceil br\rceil}\sum_{j=f_{r}^{(l)}(i)}^{f_{r}^{(u)}(i)}\sum_{j'=f_{r}^{(l)}(i')}^{f_{r}^{(u)}(i')}\int_{\left[0,1\right]^{2}}g(x+i',y+j')g(i,j)\\
&\times H_{\kappa}(\xi(x+i',y+j'))H_{\kappa}(\xi(i,j))dydx,
\end{align*}
and
\begin{align*}
D_{r}^{(3)}&=\mathbb{E}\sum_{i,i'=\lfloor ar\rfloor}^{\lceil br\rceil}\sum_{j=f_{r}^{(l)}(i)}^{f_{r}^{(u)}(i)}\sum_{j'=f_{r}^{(l)}(i')}^{f_{r}^{(u)}(i')}g(i,j)g(i',j')H_{\kappa}(\xi(i,j))H_{\kappa}(\xi(i',j')).
\end{align*}   	

Now, using (\ref{eq1}) and Assumption 1 we can rewrite the first term in (\ref{eq13}) as follows
\begin{align*}    
D_{r}^{(1)}&=\kappa!\sum_{i,i'=\lfloor ar\rfloor}^{\lceil br\rceil}\sum_{j=f_{r}^{(l)}(i)}^{f_{r}^{(u)}(i)}\sum_{j'=f_{r}^{(l)}(i')}^{f_{r}^{(u)}(i')}\int_{\left[0,1\right]^{4}}g(x+i,y+j)g(x'+i',y'+j')\\
&\times B^{\kappa}(\Vert (x'+i'-(x+i),y'+j'-(y+j))\Vert)dydy'dxdx'\\
&=\kappa!\int_{A(r)}\int_{A(r)}\dfrac{g(x,y)g(x',y')L^{\kappa}(\Vert (x'-x,y'-y)\Vert)dydy'dxdx'}{\left((x'-x)^{2}+(y'-y)^{2}\right)^{\alpha \kappa/2}}.
\end{align*}

Using change of variables (\ref{eq15}) and elementary computations, we obtain
\begin{align*}
D_{r}^{(1)}&= \kappa!r^{4-\kappa\alpha}g^{2}(r,r)
\int_{r^{-1}A(r)}\int_{r^{-1}A(r)}\dfrac{g(ru_{1},ru_{2})g(rv_{1},rv_{2})}{g^{2}(r,r)}\\
&\times \dfrac{L^{\kappa}(r\Vert (u_{1}-v_{1},u_{2}-v_{2})\Vert)du_{1}dv_{1}du_{2}dv_{2}}{((v_{1}-u_{1})^{2}+(v_{2}-u_{2})^{2})^{\alpha \kappa/2}}.
\end{align*}
Adding and subtracting $g^{*}(u_{1},u_{2})$ and $g^{*}(v_{1},v_{2})$ inside the integrals, we obtain
\begin{align}\label{eq16}
D_{r}^{(1)}= \kappa!r^{4-\kappa\alpha}g^{2}(r,r)(I_{1}+2I_{2}+I_{3}),
\end{align}
where
\begin{align*}
I_{1}&=\int_{r^{-1}A(r)}\int_{r^{-1}A(r)}\left[\dfrac{g(ru_{1},ru_{2})}{g(r,r)}-g^{*}(u_{1},u_{2})\right]\left[ \dfrac{g(rv_{1},rv_{2})}{g(r,r)}-g^{*}(v_{1},v_{2})\right]\\
&\times\frac{L^{\kappa}(r\Vert (u_{1}-v_{1},u_{2}-v_{2})\Vert)du_{1}dv_{1}du_{2}dv_{2}}{((v_{1}-u_{1})^{2}+(v_{2}-u_{2})^{2})^{\alpha \kappa/2}},
\end{align*}
\begin{align*}
I_{2}&=\int_{r^{-1}A(r)}\int_{r^{-1}A(r)}\left[\dfrac{g(ru_{1},ru_{2})}{g(r,r)}-g^{*}(u_{1},u_{2})\right]g^{*}(v_{1},v_{2})\\
&\times\frac{L^{\kappa}(r\Vert (u_{1}-v_{1},u_{2}-v_{2})\Vert)du_{1}dv_{1}du_{2}dv_{2}}{((v_{1}-u_{1})^{2}+(v_{2}-u_{2})^{2})^{\alpha \kappa/2}},\\
\text{and}\\
I_{3}&=\int_{r^{-1}A(r)}\int_{r^{-1}A(r)}\dfrac{L^{\kappa}(r\Vert (u_{1}-v_{1},u_{2}-v_{2})\Vert)g^{*}(u_{1},u_{2})g^{*}(v_{1},v_{2})du_{1}dv_{1}du_{2}dv_{2}}{((v_{1}-u_{1})^{2}+(v_{2}-u_{2})^{2})^{\alpha \kappa/2}}.
\end{align*}

Let us analyse each term $I_{i},\ i=1,2,3$, separately. The term $I_{1}$ can be estimated as
\begin{align*}
I_{1}&\leq\int_{r^{-1}A(r)}\int_{r^{-1}A(r)}\bigg\vert
\dfrac{g(ru_{1},ru_{2})}{g(r,r)}-g^{*}(u_{1},u_{2})\bigg\vert\bigg\vert \dfrac{g(rv_{1},rv_{2})}{g(r,r)}-g^{*}(v_{1},v_{2})\bigg\vert\\
&\times
\frac{L^{\kappa}(r\Vert (u_{1}-v_{1},u_{2}-v_{2})\Vert)du_{1}dv_{1}du_{2}dv_{2}}{((v_{1}-u_{1})^{2}+(v_{2}-u_{2})^{2})^{\alpha \kappa/2}}.
\end{align*}
As for each $\varepsilon$ there exists $r_{0}>0$ such that $r^{-1}A(r)\subseteq r^{-1}\Delta_{2}(r(1+\varepsilon))=\Delta_{2\varepsilon}$ for all $r>r_{0}$, then
\begin{align*}
I_{1}&\leq\int_{\Delta_{2\varepsilon}}\int_{\Delta_{2\varepsilon}}\bigg\vert
\dfrac{g(ru_{1},ru_{2})}{g(r,r)}-g^{*}(u_{1},u_{2})\bigg\vert\bigg\vert\dfrac{g(rv_{1},rv_{2})}{g(r,r)}-g^{*}(v_{1},v_{2})\bigg\vert\\
&\times\frac{L^{\kappa}(r\Vert (u_{1}-v_{1},u_{2}-v_{2})\Vert)du_{1}dv_{1}du_{2}dv_{2}}{((v_{1}-u_{1})^{2}+(v_{2}-u_{2})^{2})^{\alpha \kappa/2}}\\
&\leq\sup_{\substack{(u_{1}, u_{2})\in\Delta_{2\varepsilon}\\}}
\bigg\vert\dfrac{g(ru_{1},ru_{2})}{g(r,r)}-g^{*}(u_{1},u_{2})\bigg\vert^{2}
\int_{\Delta_{2\varepsilon}}\int_{\Delta_{2\varepsilon}}\frac{L^{\kappa}(r\Vert (u_{1}-v_{1},u_{2}-v_{2})\Vert)du_{1}dv_{1}du_{2}dv_{2}}{((v_{1}-u_{1})^{2}+(v_{2}-u_{2})^{2})^{\alpha \kappa/2}}.
\end{align*}

To estimate the above integral, we consider the uniform distribution on $\Delta_{2\varepsilon}$ with probability density function $\vert\Delta_{2\varepsilon}\vert^{-1}\chi_{\Delta_{2\varepsilon}}(x)$, $x\in\mathbb{R}^{2}$, where $\chi_{E}(\cdot)$ is the indicator function of a set $E\subseteq\mathbb{R}^{2}$. Let $U=(u_{1},u_{2})$ and $V=(v_{1},v_{2})$ be two random points which are independent and uniformly distributed inside the set $\Delta_{2\varepsilon}$. We denote by $\psi_{\Delta_{2\varepsilon}}(\rho),\ \rho\geq 0$, the pdf of the distance $\rho=\Vert U-V\Vert$. Note that in this case $\psi_{\Delta_{2\varepsilon}}(\rho)=0$ if $\rho>diam(\Delta_{2\varepsilon})$, and the Jacobian is equal to $\vert{J}\vert=C\rho,\ C>0$. Hence, for $\alpha_{0}<1:$
\[
\int_{\Delta_{2\varepsilon}}\int_{\Delta_{2\varepsilon}}\dfrac{du_{1}dv_{1}du_{2}dv_{2}}{((v_{1}-u_{1})^{2}+(v_{2}-u_{2})^{2})^{\alpha_{0}}}
=\vert\Delta_{2\varepsilon}\vert^{2}E\left(\chi\left(\Vert x-y\Vert\leq diam(\Delta_{2\varepsilon})\right)\Vert x-y\Vert^{-2\alpha_{0}}\right)
\]
\[
\qquad\qquad\qquad\qquad\qquad=C\vert\Delta_{2\varepsilon}\vert^{2}\notag\int_{0}^{diam(\Delta_{2\varepsilon})}\rho^{1-2\alpha_{0}}\psi_{\Delta_{2\varepsilon}}(\rho)d\rho
\]
\[
\qquad\qquad\leq C\vert\Delta_{2\varepsilon}\vert \int_{0}^{diam(\Delta_{2\varepsilon})}\rho^{1-2\alpha_{0}}d\rho=C\vert\Delta_{2\varepsilon}\vert \dfrac{(diam(\Delta_{2\varepsilon}))^{2-2\alpha_{0}}}{2-2\alpha_{0}}.
\]

As $\vert\Delta_{2\varepsilon}\vert=\left(1+\varepsilon\right)^{2}\vert\Delta_{2}\vert$ and $diam\left(\Delta_{2\varepsilon}\right)=\left(1+\varepsilon\right)diam\left(\Delta_{2}\right)$
\begin{align}\label{eq17}
&\int_{\Delta_{2\varepsilon}}\int_{\Delta_{2\varepsilon}}\dfrac{du_{1}dv_{1}du_{2}dv_{2}}{((v_{1}-u_{1})^{2}+(v_{2}-u_{2})^{2})^{\alpha_{0}}}\leq C\left(1+\varepsilon\right)^{4-2\alpha_{0}}\vert\Delta_{2}\vert
\dfrac{(diam(\Delta_{2}))^{2-2\alpha_{0}}}{2-2\alpha_{0}}<\infty.
\end{align}
It follows from Assumption \ref{ass1} that $L(\cdot)$ is locally bounded and by Theorem 1.5.3 in~\cite{bingham1989regular} for an arbitrary $\delta>0$ there exists $r_{0}$ and $C>0$ such that for all $r>r_{0}$
\begin{align*}
\sup_{\substack{0\leq s\leq r\cdot diam(\Delta_{2\varepsilon})\\}}\frac{s^{\delta}L(s)}{(r\cdot diam(\Delta_{2\varepsilon}))^{\delta} L(r\cdot diam(\Delta_{2\varepsilon}))}\leq C.
\end{align*}
Therefore, for all $r>r_{0}$
\begin{align*}
&\int_{\Delta_{2\varepsilon}}\int_{\Delta_{2\varepsilon}}\frac{L^{\kappa}(r\Vert (u_{1}-v_{1},u_{2}-v_{2})\Vert)du_{1}dv_{1}du_{2}dv_{2}}{((v_{1}-u_{1})^{2}+(v_{2}-u_{2})^{2})^{\alpha \kappa/2}}=\int_{\Delta_{2\varepsilon}}\int_{\Delta_{2\varepsilon}}\dfrac{(r\Vert (u_{1}-v_{1},u_{2}-v_{2})\Vert)^{\delta}}{r^{\delta}}\\
&\times\frac{L^{\kappa}(r\Vert (u_{1}-v_{1},u_{2}-v_{2})\Vert)du_{1}dv_{1}du_{2}dv_{2}}{((v_{1}-u_{1})^{2}+(v_{2}-u_{2})^{2})^{\frac{\delta+\kappa\alpha}{2}}}\leq \int_{\Delta_{2\varepsilon}}\int_{\Delta_{2\varepsilon}}\dfrac{du_{1}du_{2}dv_{1}dv_{2}}{((v_{1}-u_{1})^{2}+(v_{2}-u_{2})^{2})^{\frac{\delta+\kappa\alpha}{2}}}\\
&\times\sup_{\substack{(u_{1},u_{2})\in\Delta_{2\varepsilon}\\(v_{1},v_{2})\in\Delta_{2\varepsilon}}}\dfrac{L^{\kappa}(r\Vert (u_{1}-v_{1},u_{2}-v_{2})\Vert)(r\Vert (u_{1}-v_{1},u_{2}-v_{2})\Vert)^{\delta}}{r^{\delta}}\\
&\leq\sup_{\substack{0\leq s\leq r\cdot diam(\Delta_{2\varepsilon})}}\dfrac{L^{\kappa}(s)s^{\delta}(diam(\Delta_{2\varepsilon}))^{\delta}L^{\kappa}(r\cdot diam(\Delta_{2\varepsilon}))}{(r\cdot diam(\Delta_{2\varepsilon}))^{\delta}L^{\kappa}(r\cdot diam(\Delta_{2\varepsilon}))}\\
&\times\int_{\Delta_{2\varepsilon}}\int_{\Delta_{2\varepsilon}}\dfrac{du_{1}du_{2}dv_{1}dv_{2}}{((v_{1}-u_{1})^{2}+(v_{2}-u_{2})^{2})^{\frac{\delta+\kappa\alpha}{2}}}\\
&\leq\bigg(\sup_{\substack{0\leq s\leq r\cdot diam(\Delta_{2\varepsilon})}}\dfrac{L(s)s^{\delta/\kappa}}{(r\cdot diam(\Delta_{2\varepsilon}))^{\delta/\kappa}L(r\cdot diam(\Delta_{2\varepsilon}))}\bigg)^{\kappa}\\
&\times\left(diam(\Delta_{2\varepsilon})\right)^{\delta}L^{\kappa}(r\cdot diam(\Delta_{2\varepsilon}))\int_{\Delta_{2\varepsilon}}\int_{\Delta_{2\varepsilon}}\dfrac{du_{1}du_{2}dv_{1}dv_{2}}{((v_{1}-u_{1})^{2}+(v_{2}-u_{2})^{2})^{\frac{\delta+\kappa\alpha}{2}}}.
\end{align*}
As $\lim_{\substack{r\to\infty}}L(rdiam(\Delta_{2\varepsilon}))/L(r)=1$, one obtains for sufficiently large $r$ that
\begin{align*}
\int_{\Delta_{2\varepsilon}}\int_{\Delta_{2\varepsilon}}\frac{L^{\kappa}(r\Vert (u_{1}-v_{1},u_{2}-v_{2})\Vert)du_{1}dv_{1}du_{2}dv_{2}}{((v_{1}-u_{1})^{2}+(v_{2}-u_{2})^{2})^{\alpha \kappa/2}}
\leq \int_{\Delta_{2\varepsilon}}\int_{\Delta_{2\varepsilon}}\dfrac{CL^{\kappa}(r)du_{1}dv_{1}du_{2}dv_{2}}{((v_{1}-u_{1})^{2}+(v_{2}-u_{2})^{2})^{\frac{\delta+\kappa\alpha}{2}}}.
\end{align*}

It follows from the condition $\alpha\kappa<2$ that there exits $\delta>0$ such that $\frac{\delta+\kappa\alpha}{2}<1$. Then, applying the upper bound in~(\ref{eq17}) to the right hand side of the inequality and selecting $\delta$ such that $\frac{\delta+\kappa\alpha}{2}<1$ we obtain for sufficiently large $r_{0}$ that
\begin{align}\label{eq16new}
\int_{\Delta_{2\varepsilon}}\int_{\Delta_{2\varepsilon}}\frac{L^{\kappa}(r\Vert (u_{1}-v_{1},u_{2}-v_{2})\Vert)du_{1}dv_{1}du_{2}dv_{2}}{((v_{1}-u_{1})^{2}+(v_{2}-u_{2})^{2})^{\alpha \kappa/2}}\leq CL^{\kappa}(r),\quad r>r_{0}.
\end{align}

By Assumption~\ref{ass4}, we get $$I_{1}\leq CL^{\kappa}(r)o(1),\quad r\to\infty.$$

Similarly, by Remark~\ref{rem1} we get
\begin{align*}
I_{2}&\leq \int_{r^{-1}A(r)}\int_{r^{-1}A(r)}\bigg\vert\dfrac{g(ru_{1},ru_{2})}{g(r,r)}-g^{*}(u_{1},u_{2})\bigg\vert\vert g^{*}(v_{1},v_{2})\vert\\
&\times\dfrac{L^{\kappa}(r\Vert (u_{1}-v_{1},u_{2}-v_{2})\Vert)du_{1}dv_{1}du_{2}dv_{2}}{((v_{1}-u_{1})^{2}+(v_{2}-u_{2})^{2})^{\alpha \kappa/2}}\\
&\leq\int_{\Delta_{2\varepsilon}}\int_{\Delta_{2\varepsilon}}\bigg\vert\dfrac{g(ru_{1},ru_{2})}{g(r,r)}-g^{*}(u_{1},u_{2})\bigg\vert \vert g^{*}(v_{1},v_{2})\vert\\
&\times\dfrac{L^{\kappa}(r\Vert (u_{1}-v_{1},u_{2}-v_{2})\Vert)du_{1}dv_{1}du_{2}dv_{2}}{((v_{1}-u_{1})^{2}+(v_{2}-u_{2})^{2})^{\alpha \kappa/2}}\leq\sup_{\substack{(u_{1},u_{2})\in\Delta_{2\varepsilon}}}
\bigg\vert\dfrac{g(ru_{1},ru_{2})}{g(r,r)}-g^{*}(u_{1},u_{2})\bigg\vert\\
&\times\int_{\Delta_{2\varepsilon}}\int_{\Delta_{2\varepsilon}}\frac{\vert g^{*}(v_{1},v_{2})\vert L^{\kappa}(r\Vert (u_{1}-v_{1},u_{2}-v_{2})\Vert)du_{1}dv_{1}du_{2}dv_{2}}{((v_{1}-u_{1})^{2}+(v_{2}-u_{2})^{2})^{\alpha \kappa/2}}.
\end{align*}

As $g^{*}(\cdot,\cdot)$ is a bounded function  on $\Delta_{2\varepsilon}$, by the same reasons as for $I_{1}$ we obtain that $I_{2}\leq CL^{\kappa}(r)o(1),$ when $r\to\infty$.

Note, that as for each $\varepsilon$ there exists $r_{0}>0$ such that $\Delta_{2}(r)\subseteq A(r)\subseteq\Delta_{2\varepsilon}(r)$ for all $r>r_{0}$ and therefore $\Delta_{2}\subseteq r^{-1}A(r)\subseteq \Delta_{2\varepsilon}$. So, $r^{-1}A(r)$ converges to $\Delta_{2}$ when $r\to\infty$. 

Thus, for sufficiently large $r$
\[
I_{3}\sim\int_{\Delta_{2}}\int_{\Delta_{2}}\dfrac{L^{\kappa}(r\Vert (u_{1}-v_{1},u_{2}-v_{2})\Vert)g^{*}(u_{1},u_{2})g^{*}(v_{1},v_{2})du_{1}dv_{1}du_{2}dv_{2}}{((v_{1}-u_{1})^{2}+(v_{2}-u_{2})^{2})^{\alpha \kappa/2}}.
\]

Note, that $\int_{\Delta_{2}}\int_{\Delta_{2}}\dfrac{\abs{g^{*}(u_{1},u_{2})g^{*}(v_{1},v_{2})}du_{1}dv_{1}du_{2}dv_{2}}{((v_{1}-u_{1})^{2}+(v_{2}-u_{2})^{2})^{\frac{\delta+\alpha \kappa}{2}}}<\infty$. Hence, analogously to Proposition~4.1.2 in~\cite{bingham1989regular} we obtain that
\begin{align*}
I_{3}&\sim L^{\kappa}(r)\int_{\Delta_{2}}\int_{\Delta_{2}}\dfrac{g^{*}(u_{1},u_{2})g^{*}(v_{1},v_{2})du_{1}dv_{1}du_{2}dv_{2}}{((v_{1}-u_{1})^{2}+(v_{2}-u_{2})^{2})^{\alpha \kappa/2}}\sim l_{\Delta_{2}}L^{\kappa}(r),\quad r\to\infty,
\end{align*}
where
\begin{align*}
l_{\Delta_{2}}:=\int_{\Delta_{2}}\int_{\Delta_{2}}\dfrac{g^{*}(u_{1},u_{2})g^{*}(v_{1},v_{2})du_{1}dv_{1}du_{2}dv_{2}}{((v_{1}-u_{1})^{2}+(v_{2}-u_{2})^{2})^{\alpha \kappa/2}}.
\end{align*}
By combining the above results for (\ref{eq16}), we obtain
\begin{align*}
D_{r}^{(1)}\sim \kappa!r^{4-\kappa\alpha}g^{2}(r,r)L^{\kappa}(r)(l_{\Delta_{2}}+o(1)),\quad r\to\infty.
\end{align*}
Now, let us consider the second term $D_{r}^{(2)}$:
\begin{align*}    
D_{r}^{(2)}&=-2\kappa!\sum_{i=\lfloor ar\rfloor}^{\lceil br\rceil}\sum_{j=f_{r}^{(l)}(i)}^{f_{r}^{(u)}(i)}\int_{A(r)}\dfrac{g(x,y)g(i,j)L^{\kappa}(\Vert(i-x,j-y)\Vert)dxdy}{((i-x)^{2}+(j-y)^{2})^{\kappa\alpha/2}}\\
&=-2\kappa!\int_{A(r)}\sum_{i=\lfloor ar\rfloor}^{\lceil br\rceil}\sum_{j=f_{r}^{(l)}(i)}^{f_{r}^{(u)}(i)}\dfrac{g(x,y)g(i,j)L^{\kappa}(\Vert(i-x,j-y)\Vert)dxdy}{((i-x)^{2}+(j-y)^{2})^{\kappa\alpha/2}}.
\end{align*} 

Multiplying and dividing by $g^{2}(r,r)$ and using transformation (\ref{eq15}) again, one obtains
\begin{align*}    
D_{r}^{(2)}&=-2\kappa!r^{4-\kappa\alpha/2}g^{2}(r,r)\int_{r^{-1}A(r)}\sum_{i=\lfloor ar\rfloor}^{\lceil br\rceil}\sum_{j=f_{r}^{(l)}(i)}^{f_{r}^{(u)}(i)}\dfrac{g(ru_{1},ru_{2})g(\frac{i}{r}r,\frac{j}{r}r)}{((\frac{i}{r}-u_{1})^{2}+(\frac{j}{r}-u_{2})^{2})^{\kappa\alpha/2}}\\
&\times\dfrac{L^{\kappa}(r\Vert(\frac{i}{r}-u_{1},\frac{j}{r}-u_{1})\Vert)du_{1}du_{2}}{r^{2}g^{2}(r,r)}.
\end{align*}

Adding and subtracting $g^{*}(u_{1},u_{2})$ and $g^{*}\left(\frac{i}{r},\frac{j}{r}\right)$ inside the integral, we get
\begin{align}\label{eq18}
D_{r}^{(2)}=-2 \kappa!r^{4-\kappa\alpha}g^{2}(r,r)\left(\hat{I_{1}}+\hat{I_{2}}+\hat{I_{2}}^{\prime}+\hat{I_{3}}\right),
\end{align}
where
\begin{align*}
\hat{I_{1}}&=\int_{r^{-1}A(r)}\sum_{i=\lfloor ar\rfloor}^{\lceil br\rceil}\sum_{j=f_{r}^{(l)}(i)}^{f_{r}^{(u)}(i)}\bigg[\frac{g(ru_{1},ru_{2})}{g\left(r,r\right)}-g^{*}(u_{1},u_{2})\bigg]\bigg[\frac{g\left(\frac{i}{r}r,\frac{j}{r}r\right)}{g(r,r)}-g^{*}\left(\frac{i}{r},\frac{j}{r}\right)\bigg]\\
&\times\frac{L^{\kappa}\big(r\big\Vert\big(\frac{i}{r}-u_{1},\frac{j}{r}-u_{2}\big)\big\Vert\big)du_{1}du_{2}}{r^{2}\bigg(\big(\frac{i}{r}-u_{1}\big)^{2}+\left(\frac{j}{r}-u_{2}\right)^{2}\bigg)^{\alpha \kappa/2}},\\
\hat{I_{2}}&=\int_{r^{-1}A(r)}\sum_{i=\lfloor ar\rfloor}^{\lceil br\rceil}\sum_{j=f_{r}^{(l)}(i)}^{f_{r}^{(u)}(i)}\left[\dfrac{g(ru_{1},ru_{2})}{g\left(r,r\right)}-g^{*}(u_{1},u_{2})\right]g^{*}\left(\frac{i}{r},\frac{j}{r}\right)\\
&\times\frac{L^{\kappa}\big(r\big\Vert(\frac{i}{r}-u_{1},\frac{j}{r}-u_{2})\big\Vert\big)du_{1}du_{2}}{r^{2}\left(\big(\frac{i}{r}-u_{1}\big)^{2}+\big(\frac{j}{r}-u_{2}\big)^{2}\right)^{\alpha \kappa/2}},\\
\hat{I_{2}}^{\prime}&=\int_{r^{-1}A(r)}\sum_{i=\lfloor ar\rfloor}^{\lceil br\rceil}\sum_{j=f_{r}^{(l)}(i)}^{f_{r}^{(u)}(i)}g^{*}(u_{1},u_{2})\bigg[\dfrac{g\left(\frac{i}{r}r,\frac{j}{r}r\right)}{g(r,r)}-g^{*}\left(\dfrac{i}{r},\dfrac{j}{r}\right)\bigg]\\
&\times\frac{L^{\kappa}\big (r\big\Vert(\frac{i}{r}-u_{1},\frac{j}{r}-u_{2})\big\Vert\big )du_{1}du_{2}}{r^{2}\left(\big(\frac{i}{r}-u_{1}\big)^{2}+\big(\frac{j}{r}-u_{2}\big)^{2}\right)^{\alpha \kappa/2}},
\end{align*}
and
\begin{align*}
&\hat{I_{3}}=\int_{r^{-1}A(r)}\sum_{i=\lfloor ar\rfloor}^{\lceil br\rceil}\sum_{j=f_{r}^{(l)}(i)}^{f_{r}^{(u)}(i)}g^{*}(u_{1},u_{2})g^{*}\left(\frac{i}{r},\frac{j}{r}\right)
\frac{L^{\kappa}(r\big\Vert\big(\frac{i}{r}-u_{1},\frac{j}{r}-u_{2}\big)\big\Vert)du_{1}du_{2}}{r^{2}\left(\big(\frac{i}{r}-u_{1}\big)^{2}+\left(\frac{j}{r}-u_{2}\right)^{2}\right)^{\alpha \kappa/2}}.
\end{align*}
Similarly to the upper bounds for $I_{i},\ i=1,2,3$, we can obtain estimates for the terms $\hat{I_{1}},\hat{I_{2}},\hat{I_{2}}'$ and $\hat{I_{3}}$. For example, for sufficiently large $r$,
\[
\hat{I_{1}}\leq\int_{r^{-1}A(r)}\sum_{i=\lfloor ar\rfloor}^{\lceil br\rceil}\sum_{j=f_{r}^{(l)}(i)}^{f_{r}^{(u)}(i)}\bigg\vert\dfrac{g(ru_{1},ru_{2})}{g\left(r,r\right)}-g^{*}(u_{1},u_{2})\bigg\vert\bigg\vert\dfrac{g\left(\frac{i}{r}r,\frac{j}{r}r\right)}{g(r,r)}-g^{*}\left(\frac{i}{r},\frac{j}{r}\right)\bigg\vert
\]
\[
\times\frac{L^{\kappa}\left(r\big\Vert\big(\frac{i}{r}-u_{1},\frac{j}{r}-u_{2}\big)\big\Vert\right)du_{1}du_{2}}{r^{2}\left(\big(\frac{i}{r}-u_{1}\big)^{2}+\big(\frac{j}{r}-u_{2}\big)^{2}\right)^{\alpha \kappa/2}}\leq\int_{\Delta_{2\varepsilon}}\sum_{i=\lfloor ar\rfloor}^{\lceil br\rceil}\sum_{j=f_{r}^{(l)}(i)}^{f_{r}^{(u)}(i)}\left|\dfrac{g(ru_{1},ru_{2})}{g(r,r)}-g^{*}(u_{1},u_{2})\right|
\]
\[
\times\bigg| \frac{g\left(\frac{i}{r}r,\frac{j}{r}r\right)}{g(r,r)}-g^{*}\left(\frac{i}{r},\frac{j}{r}\right)\bigg|\ \frac{L^{\kappa}\big(r\big\Vert\big(\frac{i}{r}-u_{1},\frac{j}{r}-u_{2}\big)\big\Vert\big)du_{1}du_{2}}{r^{2}\left(\big(\frac{i}{r}-u_{1}\big)^{2}+\big(\frac{j}{r}-u_{2}\big)^{2}\right)^{\alpha \kappa/2}}
\]
\[
\leq\sup_{\substack{(u_{1},u_{2})\in\Delta_{2\varepsilon}}}\bigg\vert\dfrac{g(ru_{1},ru_{2})}{g(r,r)}-g^{*}(u_{1},u_{2})\bigg\vert^{2}
\sum_{i=\lfloor ar\rfloor}^{\lceil br\rceil}\sum_{j=f_{r}^{(l)}(i)}^{f_{r}^{(u)}(i)}\int_{\Delta_{2\varepsilon}}\frac{L^{\kappa}\big(r\big\Vert\big(\frac{i}{r}-u_{1},\frac{j}{r}-u_{2}\big)\big\Vert\big)du_{1}du_{2}}{r^{2}\left(\big(\frac{i}{r}-u_{1}\big)^{2}+\big(\frac{j}{r}-u_{2}\big)^{2}\right)^{\alpha \kappa/2}}
\]
\begin{align}\label{eq19new}
&\leq {r^{-2}}\vert\Delta_{2\varepsilon}(r)\vert \sup_{\substack{(u_{1},u_{2})\in\Delta_{2\varepsilon}}}\bigg\vert\dfrac{g(ru_{1},ru_{2})}{g(r,r)}-g^{*}(u_{1},u_{2})\bigg\vert^{2}
\notag\\
&\times\max_{i\in\set{\lfloor ar\rfloor,\dots,\lceil br\rceil}}\max_{j\in\set{f_{r}^{(l)}(i),\dots,f_{r}^{(u)}(i)}}\int_{\Delta_{2\varepsilon}}\frac{L^{\kappa}\big(r\big\Vert\big(\frac{i}{r}-u_{1},\frac{j}{r}-u_{2}\big)\big\Vert\big)du_{1}du_{2}}{\left(\big(\frac{i}{r}-u_{1}\big)^{2}+\big(\frac{j}{r}-u_{2}\big)^{2}\right)^{\alpha \kappa/2}}.
\end{align}
Analogously to the upper bound in~(\ref{eq17}) one obtains
\begin{align}\label{eq19n}
&\max_{\substack{i\in\{\lfloor ar\rfloor,\dots,\lceil br\rceil\}\\j\in\{f_{r}^{(l)}(i),\dots,f_{r}^{(u)}(i)\}}}\int_{\Delta_{2\varepsilon}}\frac{du_{1}du_{2}}{\left(\big(\frac{i}{r}-u_{1}\big)^{2}+\big(\frac{j}{r}-u_{2}\big)^{2}\right)^{\alpha_{0}}}\notag\\
&\leq\max_{\substack{i\in\{\lfloor ar\rfloor,\dots,\lceil br\rceil\}\\j\in\{f_{r}^{(l)}(i),\dots,f_{r}^{(u)}(i)\}}}\int_{\Delta_{2\varepsilon}}\chi\bigg(\bigg\Vert x-\bigg(\frac{i}{r},\frac{j}{r}\bigg)\bigg\Vert\leq diam\left(\Delta_{2\varepsilon}\right)\bigg)\bigg\Vert x-\left(\frac{i}{r},\frac{j}{r}\right)\bigg\Vert^{-2\alpha_{0}}dx\notag\\
&\leq\max_{\substack{i\in\{\lfloor ar\rfloor,\dots,\lceil br\rceil\}\\j\in\{f_{r}^{(l)}(i),\dots,f_{r}^{(u)}(i)\}}}\int_{\Delta_{2\varepsilon}-\left(\frac{i}{r},\frac{j}{r}\right)}\Vert y\Vert^{-2\alpha_{0}}dy\leq\int_{\Delta_{2\varepsilon}\left(2diam\left(\Delta_{2\varepsilon}\right)\right)}\Vert y\Vert^{-2\alpha_{0}}dy<\infty.
\end{align}
Using (\ref{eq19new}), (\ref{eq19n}) and similar steps to the proof of the upper bound in~(\ref{eq16new}) for sufficiently large $r$ we get
\begin{align*}
\hat{I_{1}}&\leq C\sup_{\substack{(u_{1},u_{2})\in\Delta_{2\varepsilon}}}\bigg\vert\dfrac{g(ru_{1},ru_{2})}{g(r,r)}-g^{*}(u_{1},u_{2})\bigg\vert^{2}
\dfrac{\vert\Delta_{2\varepsilon}(r)\vert}{r^{2}}L^{\kappa}(r)\\
&\leq C\vert\Delta_{2\varepsilon}\vert\sup_{\substack{(u_{1},u_{2})\in\Delta_{2\varepsilon}}}\bigg\vert\dfrac{g(ru_{1},ru_{2})}{g(r,r)}-g^{*}(u_{1},u_{2})\bigg\vert^{2}L^{\kappa}(r).
\end{align*}
It follows from Assumption~\ref{ass4} that
$\hat{I_{1}}\leq o(1)L^{\kappa}(r),\quad r\to\infty$.

Also, we have
\[
\hat{I_{2}}\leq\int_{r^{-1}A(r)}\sum_{i=\lfloor ar\rfloor}^{\lceil br\rceil}\sum_{j=f_{r}^{(l)}(i)}^{f_{r}^{(u)}(i)}\bigg\vert\dfrac{g(ru_{1},ru_{2})}{g\left(r,r\right)}-g^{*}(u_{1},u_{2})\bigg\vert\bigg\vert g^{*}\left(\dfrac{i}{r},\dfrac{j}{r}\right)\bigg\vert
\]
\[
\times\frac{L^{\kappa}\big(r\big\Vert\big(\frac{i}{r}-u_{1},\frac{j}{r}-u_{2}\big)\big\Vert\big)du_{1}du_{2}}{r^{2}\left(\big(\dfrac{i}{r}-u_{1}\big)^{2}+\big(\frac{j}{r}-u_{2}\big)^{2}\right)^{\alpha \kappa/2}}\leq\sup_{\substack{(u_{1},u_{2})\in\Delta_{2\varepsilon}}}\bigg\vert\dfrac{g(ru_{1},ru_{2})}{g(r,r)}-g^{*}(u_{1},u_{2})\bigg\vert
\]
\[
\times\int_{\Delta_{2\varepsilon}}\sum_{i=\lfloor ar\rfloor}^{\lceil br\rceil}\sum_{j=f_{r}^{(l)}(i)}^{f_{r}^{(u)}(i)}\frac{\big\vert g^{*}\big(\frac{i}{r},\frac{j}{r}\big)\big\vert L^{\kappa}\big(r\big\Vert\big(\frac{i}{r}-u_{1},\frac{j}{r}-u_{2}\big)\big\Vert\big)du_{1}du_{2}}{r^{2}\left(\big(\frac{i}{r}-u_{1}\big)^{2}+\big(\frac{j}{r}-u_{2}\big)^{2}\right)^{\alpha \kappa/2}}.
\]

As the function $g^{*}(\cdot,\cdot)$ is bounded on $\Delta_{2\varepsilon}$, analogously to~(\ref{eq19new}) we obtain for sufficiently large $r$ that
\begin{align*}
\int_{\Delta_{2\varepsilon}}\sum_{i=\lfloor ar\rfloor}^{\lceil br\rceil}\sum_{j=f_{r}^{(l)}(i)}^{f_{r}^{(u)}(i)}\frac{\big\vert g^{*}\big(\frac{i}{r},\frac{j}{r}\big)\big\vert L^{\kappa}\big(r\big\Vert\big(\frac{i}{r}-u_{1},\frac{j}{r}-u_{2}\big)\big\Vert\big)du_{1}du_{2}}{r^{2}\left(\big(\frac{i}{r}-u_{1}\big)^{2}+\big(\frac{j}{r}-u_{2}\big)^{2}\right)^{\alpha \kappa/2}}\leq CL^{\kappa}(r).
\end{align*}

Hence, $\hat{I_{2}}\leq o(1)L^{\kappa}(r),\ r\to\infty.$
Similarly, $\hat{I_{2}}^{\prime}\leq o(1)L^{\kappa}(r),$ as $r\to\infty$. Also, analogously to the proof for $I_{3},$ one obtains $\hat{I_{3}}\sim l_{\Delta_{2}}L^{\kappa}(r),\ r\to\infty$.

By combining these results for (\ref{eq18}), we have
\begin{align*}
D_{r}^{(2)}\sim-2\kappa!r^{4-\kappa\alpha}L^{\kappa}(r)g^{2}(r,r)(l_{\Delta_{2}}+o(1)),\qquad r\to\infty.
\end{align*}
Using similar arguments as for the sums in $D^{(2)}_{r}$ we obtain
\begin{align*}
D_{r}^{(3)}\sim\kappa!r^{4-\kappa\alpha}L^{\kappa}(r)g^{2}(r,r)(l_{\Delta_{2}}+o(1)),\qquad r\to\infty.
\end{align*}

By dividing the estimates of $ D_{r}^{(i)}, \ i = 1,2,3$, by $ r^{4-\kappa\alpha} L^{\kappa}(r)g^{2}(r, r)$, for sufficiently large $r$ we get
\[
\frac{ D_{r}^{(i)}}{r^{4-\kappa\alpha} L^{\kappa}(r)g^{2}(r, r)}\sim \kappa!(l_{\Delta_{2}}+o(1)),\quad i=1,3,\quad r\to\infty,
\] 
and
\[
\frac{ D_{r}^{(2)}}{r^{4-\kappa\alpha} L^{\kappa}(r)g^{2}(r, r)}\sim -2\kappa!(l_{\Delta_{2}}+o(1))\qquad r\to\infty.
\]

Hence, the ratio in~(\ref{eq7}) converges to zero, which completes the proof.\qed
\begin{flushleft}
	\textit{Proof of Theorem~\rm{\ref{theonew7}}}. Note, that from the isonormal spectral representation~(\ref{eq0}) and the It{\^o} formula
\end{flushleft}
\begin{align*}
H_{\kappa}(\xi(x))=\int_{\mathbb{R}^{n\kappa}}^{\prime}e^{i\langle \lambda_{1}+\dots+\lambda_{\kappa},x \rangle}\prod_{j=1}^{\kappa}\sqrt{\varphi(\lambda_{j})}W(d\lambda_{j})
\end{align*}
it follows that
\begin{align*}
Z^{(c)}_{r,\kappa}=\dfrac{c_{1}^{-\kappa/2}(n,\alpha)}{r^{n-\alpha \kappa/2} g\left(r\textit{1}_{n}\right)L^{\kappa/2}(r)}\int_{\Delta_{n}(r)}\int_{\mathbb{R}^{n\kappa}}^{\prime}g(x)e^{i\langle\lambda_{1}+\dots+\lambda_{\kappa},x\rangle}\prod_{j=1}^{\kappa}\sqrt{\varphi(\Vert\lambda_{j}\Vert)}W(d\lambda_{j})dx.
\end{align*}
Using the transformation $x^{\prime}=x/r$ we get
\begin{equation}\label{eq23}
Z^{(c)}_{r,\kappa}=\dfrac{r^{n}c_{1}^{-\kappa/2}(n,\alpha)}{r^{n-\alpha \kappa/2}L^{\kappa/2}(r)}\int_{\Delta_{n}}\int_{\mathbb{R}^{n\kappa}}^{\prime}\dfrac{g(rx^{\prime})}{g\left(r\textit{1}_{n}\right)}e^{i\langle\lambda_{1}r+\dots+\lambda_{\kappa}r,x^{\prime}\rangle}\prod_{j=1}^{\kappa}\sqrt{\varphi(\Vert\lambda_{j}\Vert)}W(d\lambda_{j})dx^{\prime}.
\end{equation}

Note that, for any fixed real number $r>0$ the function $\abs{\dfrac{g(rx)}{g\left(r\textit{1}_{n}\right)}}$ is bounded on $\Delta_{n}$. Also,  for sufficiently large $r$, it follows
by Assumption~\ref{ass5} that the function $g^{*}(\cdot)$ is bounded on $\Delta_{n}$ and therefore
\begin{align*} \int_{\Delta_{n}}\abs{\frac{g(rx)}{g(r\textit{1}_{n})}}dx&=\int_{\Delta_{n}}\abs{\frac{g(rx)}{g(r\textit{1}_{n})}-g^{*}(x)+g^{*}(x)}dx\\ &\leq\vert\Delta_{n}\vert\sup_{\substack{x\in\Delta_{n}\\}} \bigg\vert\dfrac{g(rx)}{g(r\textit{1}_{n})}-g^{*}(x)\bigg\vert+\int_{\Delta_{n}}\abs{g^{*}(x)}dx=C<\infty,\quad r\to\infty. 
\end{align*}
By Assumption~\ref{ass2} it follows that $\prod_{j=1}^{\kappa}\sqrt{\varphi(\Vert\lambda_{j}\Vert)}\in L_{2}(\mathbb{R}^{n\kappa})$. So, one can apply the stochastic Fubini's theorem to interchange the order of integration in~(\ref{eq23})~(see Theorem 5.13.1 in~\cite{peccati2011wiener}), which results in
\begin{align*}
Z^{(c)}_{r,\kappa}=\dfrac{1}{r^{-\alpha \kappa/2} L^{\kappa/2}(r)c_{1}^{\kappa/2}(n,\alpha)}\int_{\mathbb{R}^{n\kappa}}^{\prime}\left(\int_{\Delta_{n}}\dfrac{g(rx)}{g\left(r\textit{1}_{n}\right)}e^{i\langle\lambda_{1}r+\dots+\lambda_{\kappa}r,x^{\prime}\rangle}dx\right)\prod_{j=1}^{\kappa}\sqrt{\varphi(\Vert\lambda_{j}\Vert)}W(d\lambda_{j}),
\end{align*}

Using the transformation $\lambda^{(j)}=r\lambda_{j},\ j=1,\dots,\kappa$, and the self-similarity of the Gaussian white noise we get
\begin{align*}
Z^{(c)}_{r,\kappa}&=\dfrac{r^{\alpha \kappa/2-n\kappa/2}}{ L^{\kappa/2}(r)c_{1}^{\kappa/2}(n,\alpha)}\int_{\mathbb{R}^{n\kappa}}^{\prime}\left(\int_{\Delta_{n}}\dfrac{g(rx)}{g\left(r\textit{1}_{n}\right)}e^{i\langle\lambda^{(1)}+\dots+\lambda^{(\kappa)},x\rangle}dx\right)\prod_{j=1}^{\kappa}\sqrt{\varphi(\Vert\lambda^{(j)}\Vert/r)}W(d\lambda^{(j)}).
\end{align*}

By Assumption~\ref{ass2} it follows that
\begin{align*}
Z^{(c)}_{r,\kappa}&=\dfrac{1}{L^{\kappa/2}(r)}\int_{\mathbb{R}^{n\kappa}}^{\prime}\left(\int_{\Delta_{n}}\dfrac{g(rx)}{g\left(r\textit{1}_{n}\right)}e^{i\langle\lambda^{(1)}+\dots+\lambda^{(\kappa)},x\rangle}dx\right)\dfrac{\prod_{j=1}^{\kappa}\sqrt{L(r/\Vert\lambda^{(j)}\Vert)}W(d\lambda^{(j)})}{\prod_{j=1}^{\kappa}\Vert\lambda^{(j)}\Vert^{(n-\alpha)/2}}\\
&=\int_{\mathbb{R}^{n\kappa}}^{\prime}\left(\int_{\Delta_{n}}\dfrac{g(rx)}{g\left(r\textit{1}_{n}\right)}e^{i\langle\lambda^{(1)}+\dots+\lambda^{(\kappa)},x\rangle}dx\right)\dfrac{\prod_{j=1}^{\kappa}\sqrt{L(r/\Vert\lambda^{(j)}\Vert)/L(r)}W(d\lambda^{(j)})}{\prod_{j=1}^{\kappa}\Vert\lambda^{(j)}\Vert^{(n-\alpha)/2}}.
\end{align*}

By the isometry property of multiple stochastic integrals
\[
R_{r}:=\mathbb{E}\big\vert Z^{(c)}_{r,\kappa}-Z_{\kappa}\big\vert^{2}\qquad\qquad\qquad\qquad\qquad\qquad\qquad\qquad\qquad\qquad
\]
\[
=\int_{\mathbb{R}^{n\kappa}}\bigg|\int_{\Delta_{n}}e^{i\langle \lambda_{1}+\dots+\lambda_{\kappa},x \rangle}
\bigg(\frac{g(rx)}{g(r\textit{1}_{n})}\prod_{j=1}^{\kappa}\sqrt{\frac{L(r/\Vert\lambda_{j}\Vert)}{L(r)}}-g^{*}(x)\bigg)dx\bigg|^{2}\prod_{j=1}^{\kappa}\Vert \lambda_{j}\Vert^{\alpha-n}\prod_{j=1}^{\kappa}d\lambda_{j}.
\]

It follows from Assumption~\ref{ass6} that $R_{r}\to 0\ \text{as}\ r\to\infty$, which completes the proof.\qed
\begin{flushleft}
	\textit{Proof of Theorem~\rm{\ref{theo7}}}. Note, that to obtain the result of the theorem it is sufficient to prove that 
\end{flushleft}
\[
\tilde{R}_{r}:=\mathbb{E}(Z^{(d)}_{r,\kappa}-Z_{\kappa})^{2}=0,\quad r\to\infty.
\]

One can estimate $\tilde{R}_{r}$ as
\begin{align*}
\tilde{R}_{r}&=\mathbb{E}(Z^{(d)}_{r,\kappa}-Z^{(c)}_{r,\kappa}+Z^{(c)}_{r,\kappa}-Z_{\kappa})^{2}\leq 2\mathbb{E}(Z^{(c)}_{r,\kappa}-Z_{\kappa})^{2}+2\mathbb{E}(Z^{(d)}_{r,\kappa}-Z^{(c)}_{r,\kappa})^{2}.
\end{align*}

By Theorem~\ref{theonew7} the term $\mathbb{E}(Z^{(c)}_{r,\kappa}-Z_{\kappa})^{2}\to 0$ as $ r\to\infty$. Also, by Theorem~\ref{theo5} the term $\mathbb{E}(Z^{(d)}_{r,\kappa}-Z^{(c)}_{r,\kappa})^{2}\to 0$ as $r\to\infty$. 
Hence, $\tilde{R}_{r}\to 0$ as $r\to\infty$, which completes the~proof.\qed
\section{Numerical Studies}\label{sec6}
This section presents numeric examples confirming that the obtained theoretical results are valid even for wider classes of cyclic long-range dependent fields with spectral singularities outside the origin. It is demonstrated that the mean square distance between additive and the corresponding integral functionals approaches zero when $r\to\infty$. We also present  simulation studies of convergence rates. All simulations were performed by using parallel computing on the NCI's high-performance computer Raijin and the R package 'RandomFileds'~(\cite{schlather2015}). A reproducible version of the code in this paper is available in the folder "Research materials" from the  website~\url{https://sites.google.com/site/olenkoandriy/}.

For numerical examples in this section we used the cyclic long-range dependent Bessel random field. Its realisations $\xi(x,y)$ on the squares $\Delta_{2}(r):=[-r,r]^{2},\ r>0$, were simulated. The covariance function of this field has the form
\[
B\big(\sqrt{x^{2}+y^{2}}\big)=2^{v}\Gamma(v+1)\frac{J_{v}\big(\sqrt{x^{2}+y^{2}}\big)}{\big(x^{2}+y^{2}\big)^{v/2}},\quad x,y\in\mathbb{R},\ 0\leq v<\frac{1}{2}.
\]
Note, that for the range $v\in[0,\frac{1}{2})$ the covariance function $B(r)$ oscilates with the amplitude $r^{-v-0.5},\ r\to\infty$, and is not integrable which means that the long-range dependent case is considered. 

As computer simulations are possible only on a discrete grid the integrals functionals in Theorem~\ref{theo4} were approximated by the Riemann sums. Each unit interval was uniformly split by equidistant points with the step length $h$. Then the corresponding approximation of $Y_{r,\kappa}^{(c)}$ is
\begin{align}\label{eqsim1}
\tilde{Y}_{r,\kappa}^{(c)}:=h^{2}\sum\limits_{i=0}^{\frac{2r}{h}}\sum\limits_{j=0}^{\frac{2r}{h}}\dfrac{g(-r+hi,-r+hj)H_{\kappa}(\xi(-r+hi,-r+hj))}{r^{2-\alpha \kappa/2}L^{\kappa/2}(r)g\left(r,r\right)}.
\end{align}

The weight function $g(x,y)=1+(x+y)^{2},\ x,\ y\in\mathbb{R}$, was used. The function $g(\cdot,\cdot)$ satisfies Assumption~\ref{ass4} and $g^{*}(x,y)=\frac{1}{4}(x+y)^{2},\ x,\ y\in\mathbb{R}$. The value $v=0$ was used to simulate Bessel random fields. In this case the Bessel covariance function oscillating and has the asymptotic hyperbolic decay rate $(\sqrt{x^{2}+y^{2}})^{-0.5}$. Hence, the long-range dependence parameter can be chosen $\alpha=0.5.$ For simulations we used $\kappa=2$ and $h=0.1$. Hence, $H_{2}(x)=x^{2}-1$ and (\ref{eqsim1}) becomes
\begin{align}\label{eqsim1new}
\tilde{Y}_{r,2}^{(c)}=(0.1)^{2}\sum\limits_{i=0}^{20r}\sum\limits_{j=0}^{20r}\dfrac{\left(1+(-2r+0.1(i+j))^{2}\right)(\xi^{2}(-r+0.1i,-r+0.1j)-1)}{r^{3/2}\left(1+4r^{2}\right)}.
\end{align}
The corresponding additive functional is 
\begin{align}\label{eqsim2new}
Y_{r,2}^{(d)}=\dfrac{\sum_{i=\lfloor -r\rfloor}^{\lceil r\rceil}\sum_{j=\lfloor -r\rfloor}^{\lceil r\rceil}\left(1+(i+j)^{2}\right)(\xi^{2}(i,j)-1)}{r^{3/2}\left(1+4r^{2}\right)}.
\end{align}
 
Each random variable (\ref{eqsim1new}) and (\ref{eqsim2new}) was simulated $ 100 $ times for $r=10,20,\dots,80$. Then, for each $r$ we calculated the sample mean square distance $\hat{L}_{2}(r)$ between $\tilde{Y}_{r,2}^{(c)}$ and $Y_{r,2}^{(d)}$.

\begin{figure}[H]
	\centering
	\includegraphics[width=0.5\linewidth,height=6.75cm]{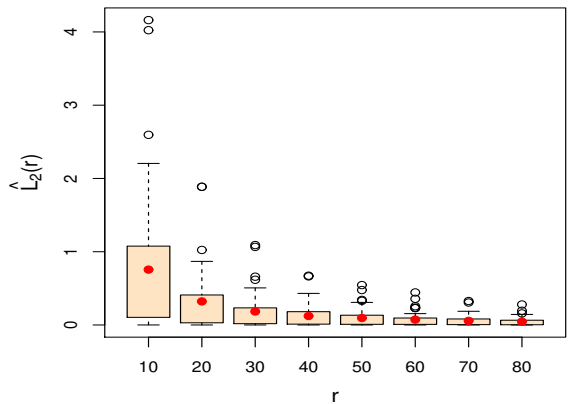}\vspace{-0.4cm}
	\caption{\label{fig3:case1} Box plots of the mean square distances between $\tilde{Y}_{r,2}^{(c)}$ and $Y_{r,2}^{(d)}$.}
\end{figure}

\begin{figure}[H]
	\centering
	\begin{subfigure}[b]{0.5\textwidth}
		\includegraphics[width=0.99\textwidth,height=6.75cm,trim={0cm 0 0 0},clip]{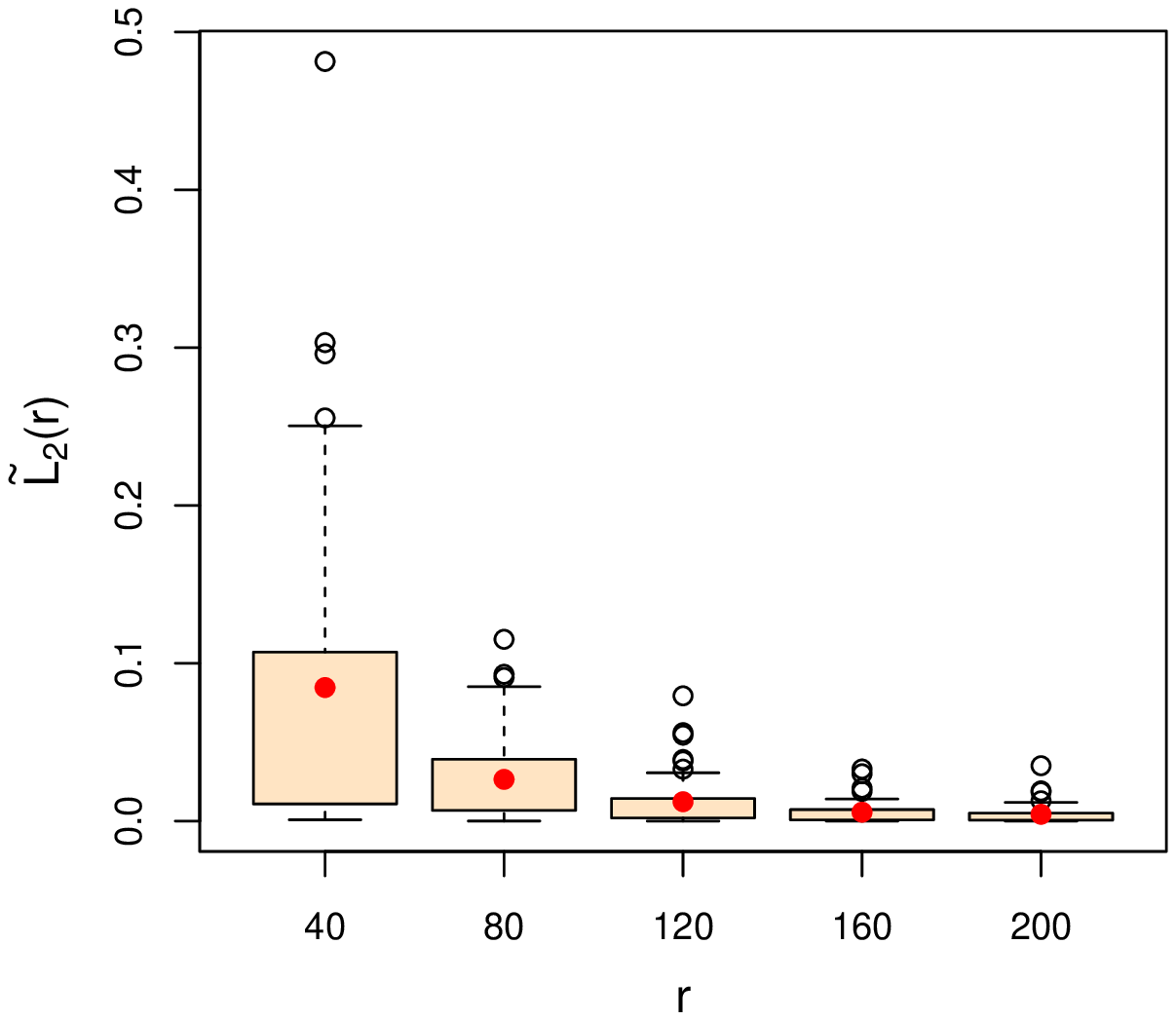}  \vspace{-0.5cm}
		\caption{}
		\label{fig4:case2a}
	\end{subfigure}
	\hfill  \hspace{-50cm}
	\begin{subfigure}[b]{0.5\textwidth} 
		\includegraphics[width=0.99\textwidth, height=6.75cm,trim={0cm 0 0 0},clip]{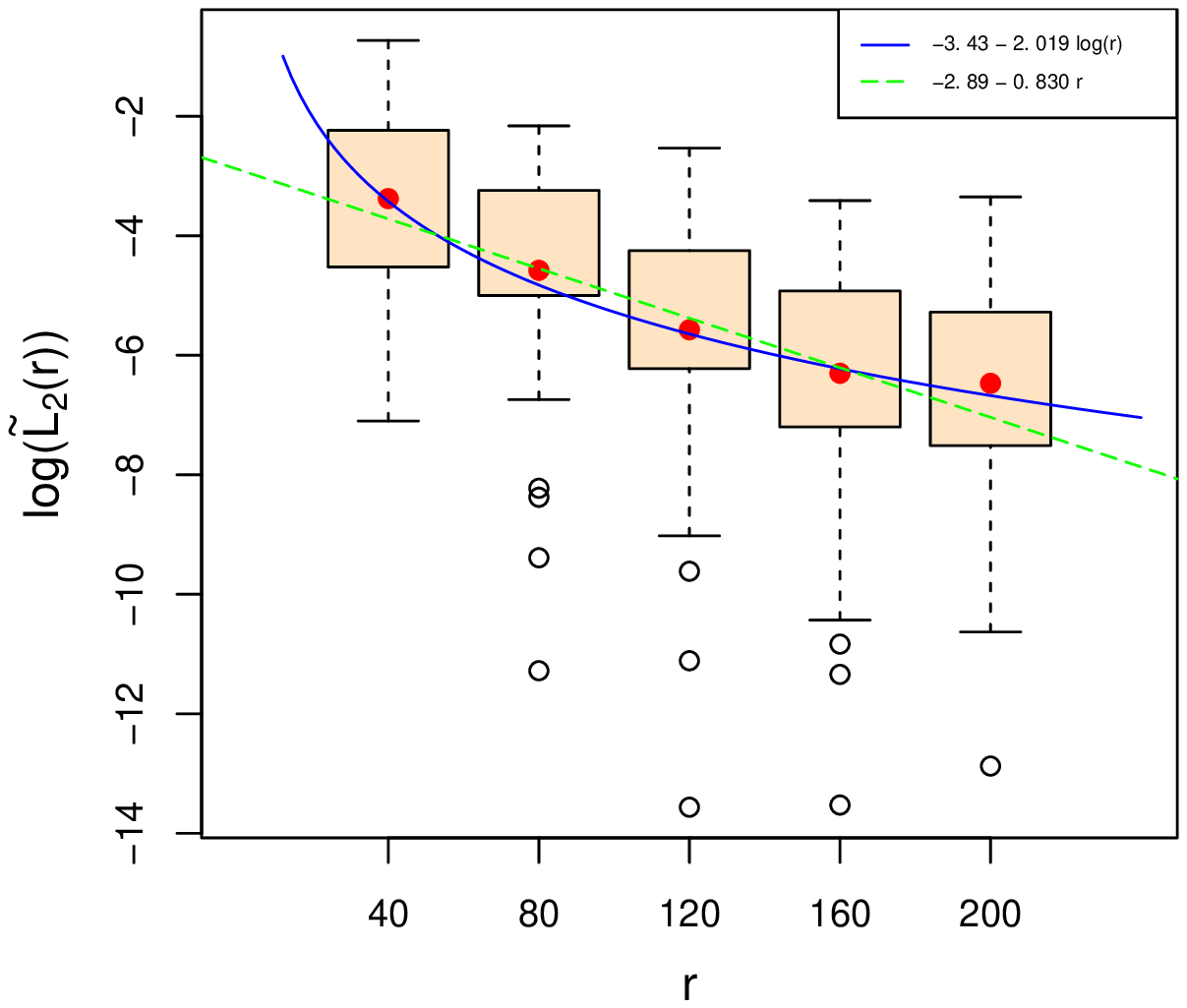} \vspace{-0.5cm}
		\caption{}
		\label{fig4:case2b}
	\end{subfigure}
	\vspace{-.2cm}
	\caption{(a) Box plots of $\tilde{L}_{2}(200,r)$, (b) Box plots of $\log(\tilde{L}_{2}(200,r))$.}\label{fig4}
\end{figure}

 Figure~\ref{fig3:case1} shows box plots of $\hat{L}_{2}(r)$ as a function of $r$ by repeating the simulation steps $50$ times, while Table~\ref{tab:Table1} shows the sample averages of $\hat{L}_{2}(r)$ for each $r$. Figure~\ref{fig3:case1} and Table~\ref{tab:Table1} confirm that the mean squared distances quickly approaches zero when $r$ increases.
 \vspace{-0.6cm} 
\begin{table}[H]
	\begin{center}
	\caption{\label{tab:Table1} Sample averages of values in the box plots from Figure~\ref{fig3:case1}.}
	\begin{tabular}{@{}lllllllll@{}}
	\toprule
	$r$ & $ 10 $ & $ 20 $ & $ 30 $ & $ 40 $ & $ 50 $ & $ 60 $ & $ 70 $ & $ 80 $\\
	\midrule
	Average of $\hat{L}_{2}(r)$ & $ 0.755 $&$ 0.322 $& $ 0.184 $&$ 0.123 $&$ 0.095 $&$ 0.073 $ & $ 0.056 $ & $ 0.046 $\\
	\bottomrule
	\end{tabular}
\end{center}
\end{table}
As the limit random variable $Z_{2}$ in Theorem~\ref{theonew7} we considered the random variable $\tilde{Y}_{R,2}^{(c)}$, where $R$ is sufficiently large. For simulations we used $\kappa=2$ and $h=0.2$. The random variables $\tilde{Y}_{R,2}^{(c)}$ and $Y_{r,2}^{(d)}$ were simulated $100$ times for $r=40,80,\dots,200$, and $R=200$. Using the simulated values, the sample mean square distance $\tilde{L}_{2}(200,r)$ between $\tilde{Y}_{200,2}^{(c)}$ and $Y_{r,2}^{(d)}$ was calculated for each~$r$. Figure~\ref{fig4} shows box plots of $\tilde{L}_{2}(200,r)$ as function of $r$ and the corresponding box plots of the logarithms of $\tilde{L}_{2}(200,r)$ by repeating the above simulation steps $50$ times. The sample averages of $\tilde{L}_{2}(200,r)$ for each $r$ is listed in Table~\ref{tab:Table3}.  Figure~\ref{fig4:case2a} and Table~\ref{tab:Table3} confirm that $\tilde{L}_{2}(200,r)$ approaches zero when $r$ increases. 

From Figure~\ref{fig4:case2b} one can see that the means form a declining slope, which suggests that the exact rate of convergence might be a power or even exponential function of $r$. By fitting the linear regression model to $\log(\tilde{L}_{2}(200,r))$ values we obtained the following $\log$-transformed models $\log(\tilde{L}_{2}(200,r))\approx -3.43 - 2.019 \log(r)$ and $\log(\tilde{L}_{2}(200,r))\approx -2.89 - 0.830r$ for the power and exponential cases respectively. The fitted models  are shown in Figure~\ref{fig4:case2b} as the solid blue (power) and dashed green (exponential) lines.
\raggedbottom
\begin{table}[H]
\begin{center}
\caption{\label{tab:Table3} Sample averages of values in the box plots from Figure~\ref{fig4:case2a}.}
\begin{tabular}{@{}llllll@{}}
		\toprule
		$r$ & $40 $ & $ 80 $ & $ 120 $ & $ 160 $ & $ 200 $\\
		\midrule
		Average of $\tilde{L}_{2}(200,r)$ & $ 0.0846 $ & $ 0.0264 $ & $ 0.0122 $ & $ 0.0056 $ & $ 0.0043 $\\
		\bottomrule
	\end{tabular}
\end{center}
\end{table}
\section{Conclusion and Directions for Future Research}\label{sec7}
This paper discussed the asymptotic behaviour of additive and integral functionals of long-range dependent random fields over increasing observation windows. It is shown that both additive and integral functionals converge to the same non-Gaussian distribution. The results were obtained under rather general assumptions on the weight functions and random fields.

The main results in Sections~\ref{sec4} and~\ref{sec5} were obtained for random fields with a spectral singularity at the zero frequency. The simulations studies in Section\ref{sec6} suggest to study the case of cyclic long-range dependent random fields that have a singularity at a non-zero frequency. 

Simulation results in Section~\ref{sec6}
suggest that the rate of convergence might be a power or exponential function of $r$. It would be interesting to obtain the exact rate of convergent for additive functionals using the approaches developed for integral functionals by~\cite{anh2017rate}.

Furthermore, the results in this paper were obtained for functionals over increasing observation windows. It would be interesting to derive similar results for high frequency asymptotics where the observation window is the same but the sampling rate increases.

Also, it would be important to obtain similar results for the case of functionals of vector data, see~\cite{omari2019reduction}.
\section*{Acknowledgements}
This research was partially supported under the Australian Research Council's Discovery Project DP160101366. This research includes extensive simulation studies using the computational cluster Raijin of the National Computational Infrastructure (NCI), which is supported by the Australian Government and La Trobe University.

\end{document}